%% file: ho.tex
\documentclass{article}

\usepackage{latexsym}
\usepackage{amssymb}
\usepackage{amsthm}
\usepackage{amsmath}

\theoremstyle{definition}

\input shorts.tex

\newcommand{\gR}{\widetilde{\mb{R}}}
\newcommand{\gC}{\widetilde{\mb{C}}}
\newcommand{\gOm}{\widetilde{\Omega}}
\newcommand{\gOmc}{\widetilde{\Omega}_\mathrm{c}}

\newcommand{\Wsc}{$\text{W}_{\text{sc}}$}
\newcommand{\WscEll}{\text{W}_{\text{sc}}\text{-Ell}}
\newcommand{\WscChar}[1]{\WscEll(#1)^{\text{c}}}
\newcommand{\xib}{\bar{\xi}}
\newcommand{\xibR}{r_\eps \xib}
\newcommand{\Pq}{\ensuremath{\widetilde{P}^2}}
\newcommand{\Pt}{\ensuremath{\widetilde{P}}}
\newcommand{\B}{\ensuremath{{\cal B}}}

\begin{document}

\title{Elliptic regularity and solvability for partial differential equations with Colombeau coefficients}
 
\author{G\"{u}nther H\"{o}rmann\footnote{Supported by FWF grant P14576-MAT,
permanent affiliation: Institut f\"ur Mathematik, Universit\"at Wien}
\ \& Michael Oberguggenberger\\
Institut f\"ur Technische Mathematik,\\
Geometrie und Bauinformatik\\
Universit\"at Innsbruck
}
\date{\today}
\maketitle

\begin{abstract}
The paper addresses questions of existence and regularity of solutions to linear
partial differential equations whose coefficients are generalized functions or
generalized constants in the sense of Colombeau. We introduce various new notions
of ellipticity and hypoellipticity, study their interrelation, and give a number
of new examples and counterexamples. Using the concept of $\G^\infty$-regularity
of generalized functions, we derive a general
global regularity result in the case of operators with constant generalized coefficients,
a more specialized result for second order operators, and a microlocal regularity
result for certain first order operators with variable generalized coefficients.
We also prove a global solvability result for operators with constant generalized
coefficients and compactly supported Colombeau generalized functions as right hand sides.
\end{abstract}



\input sec1a

\input sec2a

\input sec3a

\input sec4a

\input sec5a

\input sec6a

\input sec7b

\bibliographystyle{abbrv}
\bibliography{gueMO}

\end{document}

%% file: shorts.tex
\renewcommand{\L}{\ensuremath{\mathrm{L}}}
\newcommand{\Con}{\ensuremath{\mathrm{C}}}
\newcommand{\Cinf}{\ensuremath{\mathrm{C}^\infty}}
\newcommand{\D}{\ensuremath{{\cal D}}}
\renewcommand{\S}{\ensuremath{{\cal S}}}
\newcommand{\E}{\ensuremath{{\cal E}}}

\newcommand{\mb}[1]{\ensuremath{\mathbb{#1}}}
\newcommand{\N}{\mb{N}}

\newcommand{\R}{\mb{R}}
\newcommand{\C}{\mb{C}}

\newcommand{\G}{\ensuremath{{\cal G}}}

\newcommand{\Gt}{\ensuremath{{\cal G}_\tau}}
\newcommand{\Gc}{\ensuremath{{\cal G}_\mathrm{c}}}
\newcommand{\EM}{\ensuremath{{\cal E}_{\mathrm{M}}}}

\newcommand{\NN}{\ensuremath{{\cal N}}}
\newcommand{\g}{\mathrm{g}} 
\newcommand{\Ginf}{\ensuremath{\G^\infty}}

\newcommand{\WF}{\mathrm{WF}}
\newcommand{\singsupp}{\mathrm{sing supp}}
\newcommand{\Char}{\ensuremath{\text{Char}}}
\newcommand{\zs}{\setminus 0}
\newcommand{\CO}{\ensuremath{T^*(\Om) \zs}}

\renewcommand{\d}{\ensuremath{\partial}}
\newcommand{\diff}[1]{\frac{d}{d#1}}

\newfont{\bl}{msbm10 scaled \magstep2}

\newtheorem{thm}{Theorem}[section]
\newtheorem{lemma}[thm]{Lemma}
\newtheorem{prop}[thm]{Proposition}
\newtheorem{defn}[thm]{Definition}
\newtheorem{cor}[thm]{Corollary}
\newtheorem{rem}[thm]{Remark}
\newtheorem{ex}[thm]{Example}

\newcommand{\beq}{\begin{equation}}
\newcommand{\eeq}{\end{equation}}

\newcommand{\FT}[1]{\widehat{#1}}

\newcommand{\F}{\ensuremath{{\cal F}}}

\newcommand{\dis}[2]{\langle #1 , #2 \rangle}
\newcommand{\notmid}{\mid\kern-0.5em\not\kern0.5em}

\newcommand{\norm}[2]{{\| #1 \|}_{#2}}
\newcommand{\lone}[1]{\norm{#1}{L^1}}

\newcommand{\linf}[1]{\norm{#1}{L^\infty}}

\newcommand{\al}{\alpha}
\newcommand{\be}{\beta}
\newcommand{\ga}{\gamma}
\newcommand{\Ga}{\Gamma}
\newcommand{\de}{\delta}
\newcommand{\eps}{\varepsilon}

\newcommand{\vphi}{\varphi}

\newcommand{\la}{\lambda}

\newcommand{\om}{\omega}
\newcommand{\Om}{\Omega}
\newcommand{\sig}{\sigma}

\newcommand{\supp}{\mathop{\mathrm{supp}}}

\newcommand{\ovl}[1]{\overline{#1}}

%% file: sec1a.tex

\section{Introduction}

The purpose of this paper is to clarify a number of foundational issues in the
existence and regularity
theory for linear partial differential equations with coefficients belonging to the
Colombeau algebra of generalized functions. Questions we address are: What is a good
notion of ellipticity in the Colombeau setting? In terms of the microlocal point of view -
what should be a non-characteristic direction? What sort of regularity results are to
be expected? While regularity theory of Colombeau coefficients with classical, constant
coefficients was settled in \cite{O:92}, there is current interest in pseudodifferential
operators with Colombeau symbols \cite{Garetto:02} and microlocal analysis \cite{HdH:01,HK:01}.

Consider first a classical differential operator of order $m$ on an open set $\Om \subset \R^n$,
\begin{equation}\label{PDO}
   P(x,D) = \sum_{|\gamma| \le m} a_\gamma(x)D^\gamma
\end{equation}
where the coefficient functions $a_\gamma$ belong to $\Cinf(\Om)$ and $D^\gamma =
(- {\rm i}\,\d)^\gamma$. Such an operator is called {\em elliptic}, if its symbol
$P(x,\xi)$ satisfies an estimate of the form
\begin{equation}\label{classelliptic}
  \begin{array}{c}
   \forall K \Subset \Om\; \exists r > 0\; \exists C > 0\ {\rm such\ that}\
       \forall x \in K\; \forall \xi \in \R^n\ {\rm with}\ |\xi| \ge r: \\
       \vspace{-2mm} \\
   |P(x,\xi)| \ge C(1 + |\xi|)^m\,.
  \end{array}
\end{equation}
In the classical setting, this entails the additional property (with the same $r$
as in (\ref{classelliptic})):
\begin{equation}\label{classhypo}
  \begin{array}{c}
   \forall \alpha,\beta \in \N^n\; \exists C > 0\ {\rm such\ that}\
       \forall x \in K\; \forall \xi \in \R^n\ {\rm with}\ |\xi| \ge r: \\
       \vspace{-2mm} \\
   |\d_\xi^\alpha \d_x^\beta P(x,\xi)| \le C|P(x,\xi)|(1 + |\xi|)^{-|\alpha|}
  \end{array}
\end{equation}
which is related to the notion of hypoellipticity \cite[Section 4.1]{Hoermander:63}. Let now $u \in \G(\Om)$,
the Colombeau algebra of generalized functions on $\Om$, be a generalized solution to
the equation
\begin{equation}\label{PDE}
   P(x,D)\,u \;=\; f
\end{equation}
where $f \in \G(\Om)$. Contrary to the distributional setting, $f \in \Cinf(\Om)$ does
not entail that $u \in \Cinf(\Om)$ for elliptic equations. This follows already from the
mere fact that the ring $\gC$ of constants in $\G(\R^n)$ is strictly larger than $\C$.
For this reason, a subalgebra $\Ginf(\Om)$ of $\G(\Om)$ was introduced in \cite{O:92} with the
property that its intersection $\Ginf(\Om) \cap \D'(\Om)$ with the space of distributions
coincides with $\Cinf(\Om)$. It was shown that for constant coefficient hypoelliptic
operators and solutions $u \in \G(\Om)$ of (\ref{PDE}), $f \in \Ginf(\Om)$ implies
$u \in \Ginf(\Om)$. This regularity result easily carries over to elliptic operators
with $\Cinf$-coefficients.

In the Colombeau setting, generalized operators of the form (\ref{PDO}) with coefficients
$a_\gamma \in \G(\Om)$ arise naturally when one regularizes operators with discontinuous
coefficients or studies singular perturbations of operators with constant coefficients.
We pose the question of regularity: Under what conditions does $f \in \Ginf(\Om)$ entail
that the solution $u \in \G(\Om)$ to (\ref{PDE}) belongs to $\Ginf(\Om)$ as well?
It is obvious from the case of multiplication operators that we must require that the
coefficients in (\ref{PDO}) belong to $\Ginf(\Om)$ themselves. The Colombeau generalized
functions $a_\gamma$ are defined as equivalence classes of nets
$(a_{\gamma\eps})_{\eps \in (0,1]}$ of smooth functions. Inserting these in (\ref{PDO})
produces a representative $P_\eps(x,\xi)$ of the generalized symbol. The natural
generalization of the ellipticity condition (\ref{classelliptic}) seems to be
\begin{equation}\label{genelliptic}
  \begin{array}{c}
   \forall K \Subset \Om\; \exists N > 0\; \exists \eps_0 > 0\; \exists r > 0\;
    \exists C > 0\ {\rm such\ that}\\
    \vspace{-2mm}\\
     \forall \eps \in (0,\eps_0)\;
       \forall x \in K\; \forall \xi \in \R^n\ {\rm with}\ |\xi| \ge r: \\
       \vspace{-2mm} \\
   |P_\eps(x,\xi)| \ge C\eps^N(1 + |\xi|)^m
  \end{array}
\end{equation}
and has been proposed, among others, by \cite{NPS:98}. The starting point of this paper has been
our observation that this condition does not produce the desired elliptic regularity result.
In fact, we shall give examples of
operators $P_\eps(D)$ with constant, generalized coefficients satisfying (\ref{genelliptic})
such that the homogeneous equation (\ref{PDE}) with $f \equiv 0$ has solutions which do not
belong to $\Ginf(\Om)$. We also exhibit a multiplication operator satisfying (\ref{genelliptic}),
given by an element of $\Ginf(\R)$, whose inverse does not belong to $\Ginf(\R)$ and is not
even {\em equal in the sense of generalized distributions} (see \cite{Colombeau:84}) to an
element of $\Ginf(\R)$. The latter notion enters the picture due to results of
\cite{Garetto:02} and \cite{NPS:98}. Indeed, in \cite{Garetto:02} a strong version of
(\ref{classhypo}) for generalized symbols $P_\eps(x,\xi)$ is considered, where in addition
to property (\ref{genelliptic}) it is required (with the same $\eps_0$ and $r$) that
\begin{equation}\label{genhypo}
  \begin{array}{c}
   \forall \alpha,\beta \in \N^n\; \exists C_{\alpha\beta} > 0\ {\rm such\ that}\\
    \vspace{-2mm}\\
   \forall \eps \in (0,\eps_0)\;\forall x \in K\; \forall \xi \in \R^n\ {\rm with}\ |\xi| \ge r: \\
        \vspace{-2mm} \\
   |\d_\xi^\alpha \d_x^\beta P_\eps(x,\xi)| \le C_{\alpha\beta}|P_\eps(x,\xi)|(1 + |\xi|)^{-|\alpha|}\,.
  \end{array}
\end{equation}
As we show later, this condition does not follow from (\ref{genelliptic}). It is proved in
\cite{Garetto:02} that if a generalized operator $P_\eps(x,D)$ satisfies (\ref{genelliptic})
and (\ref{genhypo}) and $u \in \Gt(\R^n)$ is a solution to (\ref{PDE}) with
$f \in \Gt^\infty(\R^n)$, then $u$ is equal in the sense of generalized temperate distributions
to an element of $u \in \Gt^\infty(\R^n)$; here the notation $\Gt$ refers to the space
of temperate Colombeau generalized functions \cite{Colombeau:85}. The question immediately
arises whether (\ref{genelliptic}) and (\ref{genhypo}) together guarantee that $u$ actually
{\em belongs} to $\Ginf(\Om)$, if $f$ does. We give a positive answer in the case of
operators with constant generalized coefficients. Actually, a weakening of the two conditions
suffices: $r$ may depend on $\eps$ in a slowly varying fashion. In fact, various refined
notions of ellipticity are needed and will be discussed. This carries over to the
definition of a non-characteristic direction and microlocal elliptic regularity. Here
we obtain a microlocal elliptic regularity result for first order operators with
(non-constant) generalized coefficients. Contrary to the classical case, lower order
terms do matter. The general case of higher order equations with non-constant generalized
coefficients remains open and will be addressed in \cite{HOP:03}.

Having the tools for studying regularity properties at hand, we are able to answer
the question of solvability of the inhomogeneous equation with compactly supported
right hand side in the case of operators with constant generalized coefficients.
We ask for conditions on the operator $P(D)$ such that
\begin{equation}\label{solvability}
    \forall f \in \G_c(\Om): \exists u \in \G(\Om): \quad P(D) u = f\,.
\end{equation}
It was shown in \cite[Theorem 2.4]{NPS:98}
that the solvability property
(\ref{solvability}) holds for operators $P(D)$ such that
\[
   P_m(\xi_0) \text{ is invertible in } \gC \text{ for some in } \xi_0\in\R^n
\]
where $P_m(\xi)$ denotes the principal part.
Examples show that this is not the most general condition. In fact, we will prove
that solvability holds if and only if
\[
   \Pq(\xi_0) \text{ is invertible in } \gR \text{ for some in } \xi_0\in\R^n
\]
where $\Pq$ is the associated weight function. We provide an independent and short
proof based on the fundamental solution given in
\cite[Theorem 3.1.1]{Hoermander:63}.

The plan of the paper is as follows. Section 2 serves to collect the notions from
Colombeau theory needed in the sequel. In Section 3 we introduce various refinements
of the ellipticity conditions (\ref{genelliptic}) and (\ref{genhypo}) and study their
interrelations. In Section 4 we provide a number of examples and counterexamples
illustrating the situation. In Section 5 we prove the
regularity result for higher order operators with constant generalized
coefficients and give additional sufficient conditions for second order
operators. Section 6 is devoted to microlocal ellipticity properties of first order
generalized symbols and to the corresponding microlocal elliptic regularity
result, giving bounds on the wave front set of the solution.
Finally, Section 7 addresses the solvability question.

%% file: sec2a.tex
\section{Colombeau algebras}

The paper is placed in the framework of algebras of generalized functions
introduced by Colombeau in \cite{Colombeau:84, Colombeau:85}.
We shall fix the notation and discuss a number of known as well as new properties
pertinent to Colombeau generalized functions here. As a general reference we
recommend \cite{GKOS:01}.

Let $\Omega$ be an open subset of $\R^n$. The basic objects of the theory as we
use it are families
$(u_\eps)_{\eps \in (0,1]}$ of smooth functions $u_\eps \in \Cinf(\Omega)$ for
$0 < \eps \leq 1$. To
simplify the notation, we shall write $(u_\eps)_\eps$ in place of
$(u_\eps)_{\eps \in (0,1]}$ throughout.
We single out the following subalgebras:
\vspace{2mm}\\
{\em Moderate families}, denoted by $\EM(\Omega)$, are defined by the property:
\begin{equation}
  \forall K \Subset \Omega\,\forall \alpha \in \N_0^n\,
    \exists p \geq 0:\;\sup_{x\in K} |\d^\alpha u_\eps(x)|
       = O(\eps^{-p})\ \rm{as}\ \eps \to 0\,.
     \label{mofu}
\end{equation}
{\em Null families}, denoted by $\NN(\Omega)$, are defined by the property:
\begin{equation}
  \forall K \Subset \Omega\,\forall \alpha \in \N_0^n\,
     \forall q \geq 0:\;\sup_{x\in K} |\d^\alpha u_\eps(x)|
        = O(\eps^q)\ \rm{as}\ \eps \to 0\,.
      \label{nufu}
\end{equation}
In words, moderate families satisfy a locally uniform polynomial estimate as $\eps \to 0$,
together with all derivatives, while null functionals
vanish faster than any power of $\eps$ in the same situation. The null families
form a differential ideal in the collection of moderate families. The {\em Colombeau
algebra} is the factor algebra
\[
   \G(\Omega) = \EM(\Omega)/\NN(\Omega)\,.
\]
The algebra $\G(\Omega)$ just defined coincides with the {\em special Colombeau
algebra} in \cite[Definition 1.2.2]{GKOS:01}, where the notation $\G^s(\Omega)$
has been employed. However, as we will not use other variants of the algebra, we
drop the superscript $s$ in the sequel.

Restrictions of the elements of $\G(\Omega)$ to open subsets of $\Omega$ are defined
on representatives in the obvious way. One can show (see \cite[Theorem 1.2.4]{GKOS:01})
that $\Omega \to \G(\Omega)$ is a sheaf of differential algebras on $\R^n$. Thus
the support of a generalized function $u \in \G(\Omega)$ is well defined as the
complement of the largest open set on which $u$ vanishes. The subalgebra of compactly
supported Colombeau generalized functions will be denoted by $\Gc(\Omega)$.

The space of compactly supported distributions is imbedded in $\G(\Omega)$ by convolution:
\[
   \iota:\E'(\Omega) \to \G(\Omega),\;
     \iota(w) = \ \mbox{class of}\ (w \ast (\varphi_\eps)\vert_\Omega)_\eps\,,
\]
where
\begin{equation}
   \varphi_\eps(x) = \eps^{-n}\varphi\left(x/\eps\right)
               \label{molli}
\end{equation}
is obtained by scaling a fixed test function $\varphi \in \S(\R^n)$ of integral one
with all moments vanishing. By the sheaf property, this can be extended in a unique
way to an imbedding of the space of distributions $\D'(\Omega)$.

One of the main features of the Colombeau construction is the fact that this imbedding
renders $\Cinf(\Omega)$ a faithful subalgebra. In fact, given $f \in \Cinf(\Omega)$,
one can define a corresponding element of $\G(\Omega)$ by the constant imbedding
$\sigma(f) = \ \mbox{class of}\ [(\eps,x) \to f(x)]$.
Then the important equality $\iota(f) = \sigma(f)$ holds in $\G(\Omega)$. We summarize
the basic properties of the Colombeau algebra in short, referring to the literature for
details (e.g. \cite{Colombeau:85, GKOS:01, O:92, Rosinger:90}):

(a) $\G(\Omega)$ is a commutative, associative differential algebra with derivations
$\d_1, \dots, \d_n$ and multiplication $\diamond$;

(b) there is a linear imbedding of $\D'(\Omega)$ into $\G(\Omega)$;

(c) the restriction of each derivation $\d_j$ to $\D'(\Omega)$ coincides with the
usual partial derivative;

(d) the restriction of the multiplication $\diamond$ to $\Cinf(\Omega)$ coincides with
the usual product of smooth functions;

One of the achievements of the Colombeau construction is property (d): It is optimal in the
sense that in whatever algebra satisfying properties (a) - (c), the multiplication map does
not reproduce the product on $\Con^k(\Omega)$ for finite $k$, by Schwartz' impossibility result
\cite{Schwartz:54}.

We need a couple of further notions from the theory of Colombeau generalized functions.
Regularity theory in this setting is based on
the subalgebra $\G^\infty(\Omega)$ of {\em regular
generalized functions} in $\G(\Omega)$. It is defined by those elements which have
a representative satisfying
\begin{eqnarray}
  \forall K \Subset \Omega\,\exists p \geq 0\,\forall \alpha \in \N_0^n:\;
  \sup_{x\in K} |\d^\alpha u_\eps(x)| = O(\eps^{-p})\
     \ \rm{as}\ \eps \to 0\,. \label{regefu}
\end{eqnarray}
Observe the change of quantifiers with respect to formula (\ref{mofu}); locally, all derivatives
of a regular generalized function have the same order of growth in $\eps > 0$. One has
that (see \cite[Theorem 5.2]{O:92})
\[
  \G^\infty(\Omega) \cap \D'(\Omega) = \Cinf(\Omega)\,.
\]
For the purpose of describing the regularity of Colombeau generalized functions,
$\G^\infty(\Omega)$ plays the same role as $\Cinf(\Omega)$ does in the setting of distributions.

Let us also recall the {\em association relation} on the Colombeau
algebra $\G(\Omega)$. It identifies elements of $\G(\Omega)$ if they
coincide in the weak limit. That is, $u, v \in \G(\Omega)$ are called associated,
$u \approx v$, if
$ \lim_{\eps \to 0} \int\big(u_\eps(x) - v_\eps(x)\big) \psi(x)\,dx = 0 $
for all test functions $\psi \in \D(\Omega)$. This may be interpreted as the reduction of the
information on the family of regularizations to the usual distributional level.

Next, we need the notion of generalized point values. The ring of
Colombeau generalized numbers $\gC$ can be defined as
the Colombeau algebra $\G(\R^0)$, or alternatively as
the ring of constants in $\G(\R^n)$. More generally, we define
generalized points of open subsets $\Omega$ of $\R^n$ as follows:
On
\begin{equation}
  \Omega_M = \{ (x_\eps)_\eps \in \Omega^{(0,1]}: \exists p \geq 0 \ {\rm such\ that}\
  |x_\eps| = O(\eps^{-p})\ {\rm as}\ \eps \to 0\} \label{moderatepoints}
\end{equation}
we introduce an equivalence relation by
\[ (x_\eps)_\eps \sim (y_\eps)_\eps
  \ \Leftrightarrow \ \forall q \geq 0,\
  |x_\eps - y_\eps| =  O(\eps^q)\ {\rm as}\ \eps \to 0
\]
and denote by $\gOm =\Omega_M/\sim$
the set of {\em generalized points of} $\Omega$. The classes of the nets
\[
   (x_\eps)_\eps \in \gOm : \exists K\Subset \Omega \
   {\rm such\ that} \ x_\eps \in K  \ {\rm eventually \ as} \ \eps\to 0
\]
define the subset set of {\em compactly supported points}
$\gOmc$. With this notation, we clearly have that $\gC = \gR + i \gR$.

Given $u\in \G(\Omega)$ and $x \in \gOmc$, the generalized point value
$u(x) \in \gC$ is well defined as the class of $(u_\eps(x_\eps))_\eps$.
In addition, Colombeau generalized functions are characterized by their
point values (see \cite[Theorem 1.2.46]{GKOS:01}):
\[
   u=0 \mbox{ in } \G(\Omega)\ \Leftrightarrow \ u(x)= 0
   \mbox{ in } \gC \mbox{ for all } x \in \gOmc\,.
\]
The generalized numbers $\gR$ and $\gC$ form rings, but not fields. Further,
a partial order is defined on $\gR$: $r \leq s$ if there
are representatives $(r_\eps)_\eps$, $(s_\eps)_\eps$ with $r_\eps \le s_\eps$
for all $\eps$. An element $r \in \gR$ such that $0 \leq r$ with respect
to this partial order is called {\em nonnegative}. Concerning invertibility
in $\gR$ or $\gC$, we have the following results (see \cite[Theorems 1.2.38 and 1.2.39]
{GKOS:01}):

Let $r$ be an element of $\gR$ or $\gC$. Then
\begin{itemize}
\item[] $r$ is invertible if and only if
\item[] there exists some representative $(r_\eps)_\eps$
and an $m\in \N$ with $|r_\eps| \geq  \eps^m$ for sufficiently small $\eps > 0$.
\end{itemize}
Further,
\begin{itemize}
\item[] $r$ is not invertible if and only if
\item[] there exists a representative $(\tilde r_\eps)_\eps$ of $r$ and a sequence
$\eps_k\to 0$ such that $\tilde r_{\eps_k}= 0$ for all $k\in \N$, if and only if
\item[] $r$ is a zero divisor.
\end{itemize}
Concerning invertibility of Colombeau generalized function, we may state:

Let $u \in \G(\Omega)$. Then
\begin{itemize}
\item[] $u$ possesses a multiplicative inverse if and only if
\item[] there exists some representative $(u_\eps)_\eps$ such that
for every compact set $K \subset \Omega$, there is $m\in \N$ with
$\inf_{x \in K}|u_\eps(x)| \geq  \eps^m$ for sufficiently small $\eps > 0$, if and only if
\item[] $u(x)$ is invertible in $\gC$ for every $x \in \gOmc$.
\end{itemize}

We briefly touch upon the subject of linear algebra in $\gR^n$. Let $A$ be an
$(n\times n)$-matrix with coefficients in $\gR$. It defines an $\gR$-linear
map from $\gR^n$ to $\gR^n$. We have (see \cite[Lemma 1.2.41]{GKOS:01}):
\begin{itemize}
\item[] $A: \gR^n \to \gR^n$ is bijective if and only if
\item[] $\det(A)$ is an invertible element of $\gR$, if and only if
\item[] all eigenvalues of $A$ are invertible elements of $\gC$.
\end{itemize}
The last equivalence follows from the characterization of invertibility in $\gC$
above. Finally, a symmetric matrix $A$ will be called {\em positive definite},
if all its eigenvalues are nonnegative and invertible elements of $\gR$.

In order to be able to speak about symbols of differential operators,
we shall need the notion of a polynomial with generalized coefficients. The most
straightforward definition is to consider a generalized polynomial of degree $m$
as a member
\[
   \sum_{|\gamma| \leq m} a_\gamma \xi^\gamma \in \G_m[\xi]
\]
of the space of polynomials of degree $m$ in the indeterminate $\xi = (\xi_1,\dots, \xi_n)$,
with coefficients in $\G = \G(\Omega)$. Alternatively, we can and will view $\G_m[\xi]$
as the factor space
\begin{equation}\label{generalizedpolynomials}
    \G_m[\xi] = {\cal E}_{\mathrm{M},m}[\xi]/\NN_m[\xi]
\end{equation}
of families of polynomials of degree $m$ with
moderate coefficients modulo those with null coefficients. In this interpretation,
generalized polynomials $P(x,\xi)$ are represented by families
\[
   (P_\eps(x,\xi))_\eps = \Big(\sum_{|\gamma| \leq m} a_{\eps\gamma}(x) \xi^\gamma\Big)_\eps\,.
\]
Sometimes it will also be useful to regard polynomials as polynomial functions and hence
as elements of $\G(\Omega \times \R^n)$. Important special cases are
the polynomials with regular coefficients, $\G_m^\infty[\xi]$, and with
constant generalized coefficients, $\gC_m[\xi]$. The union of the spaces of
polynomials of arbitrary degree are the rings of polynomials $\G[\xi], \G^\infty[\xi]$,
and $\gC[\xi]$. Letting $D = (-i\d_1,\dots, -i\d_n)$, a differential operator $P(x,D)$
with coefficients in $\G(\Omega)$ simply is an element of $\G[D]$.

We now turn to a new notion which will be essential for the paper, the notion of
{\em slow scale nets}. Consider a moderate net of complex numbers $r = (r_\eps)_\eps
\in \gC_M$; it satisfies an estimate as exhibited in (\ref{moderatepoints}).
The {\em order of} $r$ is defined as
\[
   \kappa(r) = \sup \{q \in \R: \exists \eps_q\,\exists C_q > 0 \mbox{ such that }
                |r_\eps| \leq C_q\eps^q\,, \forall \eps \in (0,\eps_q)\}\,.
\]
\begin{lemma} \label{slowscalelemma}
Let $r = (r_\eps)_\eps \in \gC_M$. The following are equivalent:
\begin{itemize}
\item[(a)] The net $r$ has order $\kappa(r)\geq 0$\,;
\item[(b)] $\forall t \geq 0\,\exists \eps_t > 0$ such that
           $|r_\eps|^t \leq \eps^{-1}, \;\forall \eps \in (0,\eps_t)$\,;
\item[(c)] $\exists N \geq 0\,\forall t \geq 0\,\exists \eps_t > 0$ such that
           $|r_\eps|^t \leq \eps^{-N}, \;\forall \eps \in (0,\eps_t)$\,;
\item[(d)] $\exists N \geq 0\,\forall t \geq 0\,\exists \eps_t > 0\,\exists C_t > 0$ such that
           $|r_\eps|^t \leq C_t\eps^{-N}, \;\forall \eps \in (0,\eps_t)$\,;
\item[(e)] $\exists N \geq 0\,\exists \eps_0 > 0\,\forall t \geq 0\,\exists C_t > 0$ such that
           $|r_\eps|^t \leq C_t\eps^{-N}, \;\forall \eps \in (0,\eps_0)$\,;
\item[(f)] $\exists \eps_0 > 0\,\forall t \geq 0\,\exists C_t > 0$ such that
           $|r_\eps|^t \leq C_t\eps^{-1}, \;\forall \eps \in (0,\eps_0)$\,.
\end{itemize}
\end{lemma}
\begin{proof} (a) $\Leftrightarrow$ (b): If $r$ has order zero, we have that for all
$t \geq 0$ there is $\eps_t > 0$ and $C_t > 0$ such that $|r_\eps| \leq C_t\eps^{-1/2t},\,
\forall \eps \in (0, \eps_t)$. Since $C_t \leq \eps^{-1/2t}$ for sufficiently small
$\eps > 0$, assertion (b) follows. The converse direction is obvious.\\
(b) $\Rightarrow$ (c) $\Rightarrow$ (d) is clear.\\
(d) $\Rightarrow$ (b): Diminishing $\eps_t$ so that $C_t \leq \eps^{-1}$ for $\eps \in
(0,\eps_t)$ we first achieve: $\exists N \geq 0\,\forall t \geq 0\,\exists \eps_t > 0$
such that $|r_\eps|^t \leq \eps^{-(N+1)},\forall \eps \in (0,\eps_t)$. >From here we
conclude that $|r_\eps|^{t/(N+1)} \leq \eps^{-1},\forall \eps \in (0,\eps_t)$.
Since $t$ is arbitrary, the assertion (b) follows.\\
(d) $\Rightarrow$ (e): The net $r$ being moderate, there is $p\geq 0$ and $0 < \eps_0 \leq 1$
such that $|r_\eps| \leq \eps^{-p}$ for $\eps \in (0,\eps_0)$. Let $t \geq 0$. If
$\eps_t \geq \eps_0$ there is nothing to prove. Otherwise, we observe that
the net $(r_\eps)_\eps$ is bounded on $[\eps_t,\eps_0]$ by a constant
$D_t$, say. Then
\[
    |r_\eps|^t \leq \max(C_t, D_t^t)\eps^{-N},\forall \eps \in (0,\eps_0)\,.
\]
(e) $\Rightarrow$ (f) follows again by taking the $N$-th root of the inequality. Finally,
that (f) $\Rightarrow$ (d) is clear, and the proof of the lemma is complete.
\end{proof}
\begin{defn}
Nets satisfying the equivalent properties of Lemma \ref{slowscalelemma} are
termed {\em slow scale nets}.
\end{defn}
The name derives from the crucial property that
\begin{equation}
   |r_\eps|^t = O\Big(\frac{1}{\eps}\Big)\ \mbox{ as } \eps \to 0 \
      \mbox{ for all } t \geq 0\,. \label{slowscale}
\end{equation}
Another important class of nets are the {\em log-type nets}, which
are defined by the condition
\begin{equation}
   \exists \eps_0 > 0 \; \mbox{ such that } |r_\eps| \leq \log\frac{1}{\eps}\,,
      \;\forall \eps \in (0,\eps_0)\,. \label{logtype}
\end{equation}
In the same way as in Lemma \ref{slowscalelemma} one can prove that $(r_\eps)_\eps$
is of log-type if and only if:
\[
   \exists \eps_0 > 0\,\exists C > 0 \; \mbox{ such that } |r_\eps| \leq C\log\frac{1}{\eps}\,,
      \;\forall \eps \in (0,\eps_0)\,.
\]
Every log-type net is a slow scale net. In fact, every net which satisfies
\[
   \exists q \geq 0\, \exists \eps_0 > 0\; \mbox{ such that }
      |r_\eps| \leq \Big(\log\frac{1}{\eps}\Big)^q\,,
      \;\forall \eps \in (0,\eps_0)\,.
\]
is slow scale. Note, however, that $\left(\log\frac{1}{\eps}\right)^q$ does not
define a log-type net if $q > 1$, so slow scale nets need not be log-type. The
following result explores how far slow scale nets are apart from being of log-type.
\begin{lemma} \label{logtypelemma}
Let $(r_\eps)_\eps$ be a moderate net, $r_\eps \geq 0$ for all $\eps$.
The following are equivalent:
\begin{itemize}
\item[(a)] $\big(\exp(r_\eps)\big)_\eps$ is moderate;
\item[(b)] $\big(r_\eps\big)_\eps$ is log-type;
\item[(c)] $\big(r_\eps\big)_\eps$ is slow scale and:
 $\exists N \geq 0\,\exists \eps_0 > 0\,\exists C > 0$ and families
 $\big(c_{\eps t}\big)_{\eps > 0, t \in \N}$ such that
 $\sum_{t = 0}^\infty c_{\eps t} \leq 1$ and
 $r_\eps^t \leq Ct!c_{\eps t}\eps^{-N}, \;\forall \eps \in (0,\eps_0)$\,.
\end{itemize}
\end{lemma}
\begin{proof} (a)$ \Leftrightarrow$ (b): If (a) holds, there is $p \geq 0, \eps_0 > 0$
such that $\exp(r_\eps) \leq \eps^{-p}, \forall \eps \in (0,\eps_0)$. Thus
$r_\eps \leq p\log(1/\eps)$ and this means that $r_\eps \geq 0$ is log-type, as was
observed above. The converse is clear anyway.\\
(a) $\Rightarrow$ (c): By assumption, there is $\eps > 0, N \geq 0$ such that
\[
   \sum_{t = 0}^\infty \frac{r_\eps^t}{t!} \leq \eps^{-N},\ \forall \eps \in (0,\eps_0)\,.
\]
Put $c_{\eps t} = \eps^Nr_\eps^t/t!$. Then $\sum_{t = 0}^\infty c_{\eps t} \leq 1$ and
$r_\eps^t = t!c_{\eps t}\eps^{-N}, \;\forall \eps \in (0,\eps_0)$\,.\\
(c) $\Rightarrow$ (a): We have that
\[
   \exp(r_\eps) = \sum_{t = 0}^\infty \frac{r_\eps^t}{t!}
       \leq C\sum_{t = 0}^\infty c_{\eps t}\eps^{-N} \leq C\eps^{-N},\ \forall \eps \in (0,\eps_0)\,.
\]
This shows that $\exp(r_\eps)$ forms a moderate net and completes the proof.
\end{proof}
\begin{rem}
If $(r_\eps)_\eps$ is slow scale and the constants $C_t$ in Lemma \ref{slowscalelemma}(e)
satisfy $\sum_{t = 0}^\infty C_t/t! = C < \infty$, then condition (c) of
Lemma \ref{logtypelemma} is satisfied with $c_{\eps t} = C_t/(Ct!)$, hence
$(\exp(r_\eps))_\eps$ is moderate. If in addition the constants $C_t$ are actually
bounded by a constant $C'$, say, then $(r_\eps)_\eps$ is bounded. This follows
from the fact that at each fixed $\eps \in (0,\eps_0), \lim_{t \to \infty} (C'/\eps^N)^{1/t}
= 1$.
\end{rem}

%% file: sec3a.tex

\section{Notions of ellipticity in the Colombeau setting}

In this section, we study linear differential operators of the form (\ref{PDO}) with
coefficients $a_\gamma \in \Ginf(\Omega)$; throughout, $\Omega$ is an open
subset of $\R^n$ and the order of $P(x,D)$ is $m \geq 0$. The symbol and the
principal symbol, respectively,

\[
   P(x,\xi) = \sum_{|\gamma| \leq m} a_{\gamma}(x) \xi^\gamma\,, \quad
   P_m(x,\xi) = \sum_{|\gamma| = m} a_{\gamma}(x) \xi^\gamma\, \quad
\]
are elements of $\Ginf_m[\xi]$. Due to the $\Ginf$-property, every representative
$(P_\eps(x,\xi))_\eps$ satisfies an estimate from above of the following form:
\begin{equation*}
\begin{array}{c}
   \forall K \Subset \Omega\,\exists p \geq 0\,\forall \alpha,\beta \in \N_0^n\,
   \exists \eps_0 > 0\,\exists C > 0 \mbox{ such that }  \\
   \vspace{-2mm} \\
   \forall \eps \in (0,\eps_0)\,\forall x \in K\,\forall \xi \in \R^n:\\
   \vspace{-2mm} \\
   \big|\d_\xi^\alpha\d_x^\beta P_\eps(x,\xi)\big| \leq C\eps^{-p}\big(1 + |\xi|\big)^{m-|\alpha|}.
\end{array}
\end{equation*}
Various estimates from below will lead to various notions of ellipticity, which - due to
the presence of the additional parameter $\eps \in (0,1]$ - are more involved than in the
classical case. In this section, we introduce these concepts and relate them to
properties of the principal symbol. Examples distinguishing these notions and their
relevance for regularity theory will be given in the subsequent sections.
\begin{defn}
Let $P = P(x,\xi) \in \Ginf_m[\xi]$.
\begin{itemize}
\item[(a)] P is called {\em S-elliptic}, if for some representative $P_\eps(x,\xi)$ the
  following condition holds:
\begin{equation}\label{Selliptic}
  \begin{array}{c}
   \forall K \Subset \Om\; \exists N > 0\; \exists \eps_0 > 0\; \exists r > 0\ {\rm such\ that}\\
    \vspace{-2mm}\\
     \forall \eps \in (0,\eps_0)\;
       \forall x \in K\; \forall \xi \in \R^n\ {\rm with}\ |\xi| \ge r: \\
       \vspace{-2mm} \\
   |P_\eps(x,\xi)| \ge \eps^N(1 + |\xi|)^m.
  \end{array}
\end{equation}
\item[(b)] P is called {\em W-elliptic}, if
\begin{equation}\label{Welliptic}
  \begin{array}{c}
   \forall K \Subset \Om\; \exists N > 0\; \exists \eps_0 > 0\; \forall \eps \in (0,\eps_0)\;
    \exists r_\eps > 0\ {\rm such\ that}\\
    \vspace{-2mm}\\
     \forall \eps \in (0,\eps_0)\;
       \forall x \in K\; \forall \xi \in \R^n\ {\rm with}\ |\xi| \ge r_\eps: \\
       \vspace{-2mm} \\
   |P_\eps(x,\xi)| \ge \eps^N(1 + |\xi|)^m.
  \end{array}
\end{equation}
\item[(c)] P is called {\em SH-elliptic}, if it is S-elliptic and in addition
(with the same $\eps_0$ and $r$ as in (a))
\begin{equation}\label{SHelliptic}
  \begin{array}{c}
   \forall \alpha,\beta \in \N_0^n\; \exists C_{\alpha\beta} > 0\ {\rm such\ that}\\
    \vspace{-2mm}\\
   \forall \eps \in (0,\eps_0)\;\forall x \in K\; \forall \xi \in \R^n\ {\rm with}\ |\xi| \ge r: \\
        \vspace{-2mm} \\
   |\d_\xi^\alpha \d_x^\beta P_\eps(x,\xi)| \le C_{\alpha\beta}|P_\eps(x,\xi)|(1 + |\xi|)^{-|\alpha|}\,.
  \end{array}
\end{equation}
\item[(d)] P is called {\em WH-elliptic}, if it is W-elliptic and in addition
(with the same $\eps_0$ and $r_\eps$ as in (b))
\begin{equation}\label{WHelliptic}
  \begin{array}{c}
   \forall \alpha,\beta \in \N_0^n\; \exists C_{\alpha\beta} > 0\ {\rm such\ that}\\
    \vspace{-2mm}\\
   \forall \eps \in (0,\eps_0)\;\forall x \in K\; \forall \xi \in \R^n\ {\rm with}\ |\xi| \ge r_\eps: \\
        \vspace{-2mm} \\
   |\d_\xi^\alpha \d_x^\beta P_\eps(x,\xi)| \le C_{\alpha\beta}|P_\eps(x,\xi)|(1 + |\xi|)^{-|\alpha|}\,.
  \end{array}
\end{equation}
\end{itemize}
\end{defn}
\begin{rem}
(i) It is clear that if any of these conditions holds for one representative (in the
sense of (\ref{generalizedpolynomials})) of $P(x,\xi)$, then it holds for all
representatives. Indeed, if $(Q_\eps(x,\xi))_\eps$ belongs to $\NN[\xi]$, then
\begin{equation*}
\begin{array}{c}
   \forall K \Subset \Omega\,\forall q \geq 0\,\exists \eps_0 > 0 \mbox{ such that }  \\
   \vspace{-2mm} \\
   \forall \eps \in (0,\eps_0)\,\forall x \in K\,\forall \xi \in \R^n:\\
   \vspace{-2mm} \\
   \big|Q_\eps(x,\xi)\big| \leq \eps^{q}\big(1 + |\xi|\big)^{m}
\end{array}
\end{equation*}
and this entails the assertion.

(ii) The letters ``S'' and ``W'' should be reminiscent of ``strong'' and ``weak'', respectively,
while the ``H'' is intended to invoke an association with ``hypo-''. The ``weak'' conditions
differ from the ``strong'' ones by the fact that the radius $r$ from which on the basic
estimate is required to hold may grow as $\eps \to 0$.

(iii) The implications (a) $\Rightarrow$ (b), (c) $\Rightarrow$ (d) as well as
(c) $\Rightarrow$ (a) and (d) $\Rightarrow$ (b) are obvious. None of the reverse
implications hold, as will be seen by the examples in Section 4.
\end{rem}
\begin{prop}\label{principalpartprop}
Let $P(x,\xi) \in \Ginf_m[\xi]$ be an operator of order $m$. Then
$P(x,\xi)$ is W-elliptic if and only if its principal part
$P_m(x,\xi)$ is S-elliptic.
\end{prop}
\begin{proof}
Assume that $P_m$ is S-elliptic. Separating the homogeneous terms, we may write
\[
   P_\eps(x,\xi) = P_{m,\eps}(x,\xi) + P_{m-1,\eps}(x,\xi) + \dots +
   P_{0,\eps}(x,\xi)\,.
\]
By assumption and the fact that each coefficient $a_{\gamma\eps}$
is moderate, we have
\[
\begin{array}{r}
   |P_\eps(x,\xi) | \geq \eps^N(1+|\xi|)^m -
   C'\eps^{-N'}(1+|\xi|)^{m-1}\\
   \vspace{-2mm}\\
      = (1+|\xi|)^m\eps^N\big(1-C'\eps^{-N-N'}(1+|\xi|)^{-1}\big)
\end{array}
\]
for certain constants $C, N' > 0$, when $x$ varies in a relatively
compact set $K, |\xi| \geq r$ and $\eps \in (0,\eps_0)$. Defining $r_\eps$ by the property that
\[
   \big(1-C'\eps^{-N-N'}(1+|\xi|)^{-1}\big) \geq \frac{1}{2}\,,
\]
we get
\[
  |P_\eps(x,\xi) | \geq \frac{1}{2}(1+|\xi|)^m\eps^N
      \geq (1+|\xi|)^m\eps^{N+1}
\]
if $|\xi| \geq r_\eps$, as desired.

Conversely, assume that $P$ is W-elliptic. Let $\eta \in \R^n$
with $|\eta| = 1$ and choose $\xi \in \R^n$ such that $|\xi| \geq
r_\eps$ and $\eta = \xi/|\xi|$. Then
\[
   \big|P_\eps(x,\eta) + P_{m-1,\eps}(x,\eta)\frac{1}{|\xi|} + \dots +
   P_{0,\eps}(x,\eta)\frac{1}{|\xi|^m}\big| \geq \eps^N\big(1 + \frac{1}{|\xi|}\big)^m
\]
by hypothesis. By the moderation property of each term, we obtain for $x \in K$ that
\[
    P_{m-j,\eps}(x,\eta)\frac{1}{|\xi|^j} \leq \eps^{-N_j}\frac{1}{|\xi|^j} \leq \eps^{N+1}
\]
if $|\xi| \geq \max\big(r_\eps, \eps^{-N_j-N-1}\big)$. Thus
\[
   \big|P_{m,\eps}(x,\eta)\big| \geq \frac{1}{2}\eps^N
\]
for sufficiently small $\eps > 0$ and all $\eta$ with $|\eta| = 1$. We conclude that
\[
   \big|P_{m,\eps}(x,\xi)\big| \geq \frac{1}{2}\eps^N|\xi|^m \geq \eps^{N+1}(1+|\xi|)^m
\]
for $\xi \in \R^n, |\xi| \geq 1$ and sufficiently small $\eps > 0$, as required.
\end{proof}

For operators with constant (generalized) coefficients, the S-ellipticity of the principal
part can be characterized by the usual pointwise conditions, using generalized points.
\begin{prop} \label{pmdoesnotvanish}
Let $P(\xi) \in \gC_m[\xi]$ be an operator with constant coefficients in $\gC$. The following
are equivalent:
\begin{itemize}
\item[(a)] $P_m$ is S-elliptic;
\item[(b)] for each representative $P_\eps(x,\xi)$ it holds that
\begin{equation*}
  \begin{array}{c}
   \exists N > 0\; \exists \eps_0 > 0\ {\rm such\ that}\
       \forall \eps \in (0,\eps_0)\; \; \forall \xi \in \R^n: \\
       \vspace{-2mm} \\
   |P_{m,\eps}(\xi)| \ge \eps^N|\xi|^m\,;
  \end{array}
\end{equation*}
\item[(c)] $\forall \tilde{\xi} \in \gR^n: P_m(\tilde{\xi}) = 0 \Leftrightarrow \tilde{\xi} = 0$\,.
\end{itemize}
\end{prop}
\begin{proof}
(a) $\Rightarrow$ (b): If $P_m$ is S-elliptic and $\eta = \xi/|\xi|$, there is
$N \geq 0, \eps_0 > 0$ and $r > 0$ such that
\[
  |P_{m,\eps}(\eta)| \geq \eps^N\big(1 + \frac{1}{|\xi|}\big)^m
\]
for $\eps \in (0,\eps_0)$ and $|\xi| \geq r$. Letting $|\xi| \to \infty$ we obtain
$|P_{m,\eps}(\eta)| \geq \eps^N$ whenever $\eps \in (0,\eps_0)$ and $|\eta| = 1$. This in turn
implies (b).

(b) $\Rightarrow$ (a): If $|\xi| \geq 1$, we have that
\[
  |P_{m,\eps}(\xi)| \geq \eps^N|\xi|^m \geq \frac{\eps^N}{2^m}(1 + |\xi|)^m \geq \eps^{N+1}(1 + |\xi|)^m
\]
whenever $0 < \eps < \min(\eps_0,1/2^m)$.

(b) $\Rightarrow$ (c): If $\tilde{\xi} \neq 0$ in $\gR^n$ then there is $q \geq 0$ and a sequence
$\eps_k \to 0$ such that $|\xi_{\eps_k}| \geq \eps_k^q$, where $(\xi_\eps)_\eps$ is a representative
of $\tilde{\xi}$. But then (b) implies that $|P_{m,\eps_k}(\xi_{\eps_k}) \geq \eps_k^{N+mq}$ for
sufficiently large $k \in \N$, so $P_m(\tilde{\xi}) \neq 0$ in $\gC$.

(c) $\Rightarrow$ (b): The negation of (b) is:
\begin{equation*}
   \forall N > 0\; \forall \eps_0 > 0\;\exists \eta \in \R^n, |\eta| = 1\ {\rm such\ that}\
          |P_{m,\eps}(\eta)| < \eps^N.
\end{equation*}
In particular, choosing $\eps_0 = 1/N$ we obtain $\eps_N \in (0,
1/N)$ and $\eta_N$ with $|\eta_N| = 1$ such that
$|P_{m,\eps_N}(\eta_{N})| < \eps_N^N$\,. Note that $P_{m,\eps}(0)
= 0$. Define $\tilde{\eta} \in \gR^n$ as the class of
\[
    \eta_\eps = \left\{ \begin{array}{ll}
                        0\,, & {\rm if\ }\eps \not\in \{\eps_1, \eps_2, \eps_3,\dots\}\\
                        \eta_N\,, & {\rm if\ } \eps = \eps_N\ {\rm for\ some}\ N \in \N\,.
                        \end{array}\right.
\]
Then clearly $\tilde{\eta} \neq 0$ in $\gR^n$, but $P_m(\tilde{\eta}) = 0$ in $\gC$.
\end{proof}

%% file: sec4a.tex

\section{Examples}

In this section we present examples that distinguish the notions
of ellipticity defined in the previous sections and relate them to
generalized hypoellipticity. In particular, we shall scrutinize
operators $P(x,D)$ with coefficients in $\Ginf(\Omega)$ or in
$\gC$ with respect to the regularity property:
\begin{equation} \label{regproperty}
   \big[u \in \G(\Omega)\,,\;  f \in \Ginf(\Omega) \mbox{ and }
     P(x,D)u =f \mbox{ in } \G(\Omega) \big]
    \Longrightarrow u \in \Ginf(\Omega)\,.
\end{equation}
Operators that enjoy property (\ref{regproperty}) on every open subset
$\Omega \subset \R^n$ are called {\it $\Ginf$-hypoelliptic}. The examples given
here will also illuminate the range of validity of the hypoellipticity
results in the subsequent sections.
\begin{ex}\label{ex1}
The operator $P(\xi) \in \gC[\xi_1,\xi_2]$ on $\R^2$ defined by the
representative $P_\eps(\xi) = \eps\xi_1 + i \xi_2$: It is S-elliptic,
but not WH-elliptic (hence also W-elliptic, but not SH-elliptic).
Indeed, it is homogeneous of degree one and
\[
   |\eps\xi_1 + i \xi_2| = \sqrt{\eps^2\xi_1^2 + \xi_2^2} \geq
       \eps|\xi|
\]
for $\eps \in (0,1]$ and all $\xi \in \R^2$, hence $P(\xi)$ is S-elliptic
by Proposition \ref{pmdoesnotvanish}. On the other hand, the inequality
\[
    \big|\d_{\xi_2} P_\eps(\xi)\big| \leq C |\eps\xi_1 + i \xi_2|(1 + |\xi|)^{-1}
\]
entails for $\xi = (\xi_1,0)$ that
\[
    1 \leq C\eps|\xi_1|(1 + |\xi_1|)^{-1}
\]
and thus does not hold for whatever $C$ when $|\xi_1| \to \infty$. Thus there is
no family of radii $r_\eps$ which could produce the WH-ellipticity estimate
(\ref{WHelliptic}). The corresponding homogeneous differential equation on
$\Omega = \R \times (0,\infty)$,
\[
    \big(-i\eps\frac{\d}{\d x_1} + \frac{\d}{\d x_2}\big)\,u_\eps(x_1,x_2) = 0
\]
is solved by the moderate family $u_\eps(x_1,x_2) = \exp(ix_1/\eps - x_2)$ which
does not define an element of $\Ginf(\Omega)$, so $P(D)$ is not $\Ginf$-hypoelliptic.
\end{ex}
\begin{ex}\label{ex2}
The operator $P(\xi) = a^2\xi_1^2 + \xi_2^2 \in \gC[\xi_1,\xi_2]$ where $a$ is
a nonnegative element of $\gR$: We first consider the case when $a$ is not
invertible. Then $P$ is not S-elliptic. Indeed, $a$ is a zero divisor (Section 2), so
there is $b \in \gR, b \neq 0$ such that $ab = 0$. Taking $\tilde{\xi} = (b, 0) \in
\gR^2$ we have that $\tilde{\xi} \neq 0$, but $P(\tilde{\xi}) = 0$ (see Proposition
\ref{pmdoesnotvanish}). Also, if $u = u(\xi_1)$ is any element of $\G(\R^2)$ not
depending on $\xi_2$, then $P(D)\,bu = 0$, so $P(D)$ is not $\Ginf$-hypoelliptic.

Second, assume that $a$ is invertible. Then $P(\xi)$ is S-elliptic. Indeed, by the
invertibility of $a$, its representatives satisfy
an estimate from below of the form $a_\eps \geq \eps^m$ for some $m$ and sufficiently
small $\eps > 0$. Hence
\[
    |P_\eps(\xi)| \geq \min(1, a_\eps^2)(\xi_1^2 + \xi_2^2) \geq \eps^{2m}|\xi|^2\,,
\]
demonstrating the S-ellipticity by Proposition \ref{pmdoesnotvanish}.
On the other hand, if $a_\eps$ is not bounded away from zero by a positive,
real constant, then  $P(\xi)$ is not WH-elliptic. Indeed, an inequality of the
form
\[
    \big|\d_{\xi_2}^2 P_\eps(\xi)\big| \leq C |a_\eps^2\xi_1^2 + \xi_2^2|(1 + |\xi|)^{-2}
\]
means for $\xi = (\xi_1,0)$ that
\[
    2 \leq Ca_\eps\xi_1^2(1 + |\xi_1|)^{-2}\,.
\]
Letting $|\xi_1| \to \infty$ produces the lower bound $a_\eps \geq 2/C$.

Finally, the moderate family $u_\eps(x_1,x_2) = \sin(x_1/a_\eps)\sinh(x_2)$ defines
a solution of the homogeneous equation $P(D)u = 0$ in $\G(\R^2)$. When $(1/a_\eps)_\eps$
is not slow scale, $u$ does not belong to $\Ginf(\R^2)$, in which case $P(D)$
is not $\Ginf$-hypoelliptic.
\end{ex}
Ordinary differential operators may serve to illustrate a few further points:
\begin{ex} \label{ex3}
The operator $P(\xi) = - a^2 \xi^2 +1 \in \gC[\xi]$ on $\R$ where
$a$ is a nonnegative, invertible element of $\gR$: Its principal symbol
$P_m(\xi)$ is clearly S-elliptic, so $P(\xi)$ is W-elliptic. We will show that
$P(\xi)$ is not S-elliptic if the representatives $(a_\eps)_\eps$ are not bounded away
from zero by a positive, real constant. This demonstrates that the full symbol
of an operator need not inherit
the S-ellipticity property from the principal symbol (contrary to the
classical situation). Indeed, the fact that $P(1/a_\eps) = 0$ shows that an
estimate
\[
    |P_\eps(\xi)| \geq \eps^N(1 + |\xi|)^2
\]
cannot hold for $|\xi| \geq r$ with $r$ independent of $\eps$, if
$a_\eps$ has a subsequence converging to zero.

However, $P(\xi)$ is WH-elliptic, and we shall now estimate the required
radius; in fact, $r_\eps = s/a_\eps$ for arbitrary $s > 1$
does the job.
Indeed, if $r_\eps = s / a_\eps$ with $s > 1$ and $|\xi| \geq
r_\eps$, we have $|P_\eps(\xi)| = a_\eps^2 \xi^2 - 1 > 0$ and
\[
  \big|\frac{\d^2}{\d\xi^2} P_\eps(\xi)\big| = 2a_\eps^2
      = 2  \frac{|P_\eps(\xi)|}{\xi^2 - 1/a_\eps^2 }
      \leq 2  \frac{|P_\eps(\xi)|}{\xi^2(1 - 1/s^2)}
      \leq  c |P_\eps(\xi)|( 1 + |\xi|)^{-2}
\]
where $c = 8s^2/(s^2 - 1)$ and $|\xi| \geq \max(r_\eps, 1)$.
A similar estimate holds for the first derivative $\d_\xi P_\eps(\xi)$,
thus the second condition in the definition of WH-ellipticity holds. What
concerns the first, we must come up with $N \geq 0$ such that
\begin{equation} \label{ex_W_ell}
  a_\eps^2 \xi^2 - 1 \geq \eps^N (1 + |\xi|)^2
\end{equation}
for $\eps\in (0,\eps_0)$ and $|\xi| \geq r_\eps$. If we prove that
\[
  a_\eps^2 \xi^2 - 1 \geq \eps^N \xi^2
\]
in this range, then (\ref{ex_W_ell}) follows by possibly enlarging $N$.
Since $a$ is invertible, there is $q \geq 0$ such that
$a_\eps \geq \eps^q$ for sufficiently small $\eps$. We may choose $N =
2q + 1$, as is seen from the calculation
\[
  a_\eps^2 \xi^2 - \eps^N \xi^2 \geq
     (a_\eps^2 - \eps^{2q + 1}) \frac{s^2}{a_\eps^2}
  = s^2 (1 - \eps (\frac{\eps^q}{a_\eps})^2) \geq s^2 (1 - \eps)
     \geq 1
\]
for sufficiently small $\eps$ and $|\xi| \geq r_\eps$.

The corresponding homogeneous equation
\[
  (a_\eps^2 \frac{d^2}{dx^2} + 1) u_\eps = 0
\]
has the solution $u_\eps(x) = \sin(x / a_\eps)$ which does not define an element
of $\Ginf(\Om)$, unless $1/a_\eps$ is slow scale. This shows that the
WH-property alone does not guarantee $\Ginf$-hypoellipticity of the operator
and suggests that conditions on the radius $r_\eps$ have to enter
the picture (as will be expounded in Section 5).
\end{ex}
\begin{rem} \label{remarksinxovereps}
The case $a_\eps = \eps$ deserves some more attention. In this case the class
$u \in \G(\R)$ of $u_\eps(x) = \sin(x/\eps)$, while being a non-regular solution to the
homogeneous equation $a^2 d^2 u(x)/dx^2 +1 = 0$, is actually equal to zero in the sense of
generalized distributions, that is,
\[
   \int u(x)\psi(x)\,dx = 0 \mbox{ in } \gC
\]
for all $\psi \in \D(\R)$. Indeed,
\[
   \big|\int \sin(\frac{x}{\eps})\psi(x)\,dx\big|
   = \big|\eps^{4q}\int \sin(\frac{x}{\eps})\psi^{(4q)}(x)\,dx\big|
   = O(\eps^{4q})
\]
for every $q \geq 0$. This shows, in particular, that an element
of $\G(\R)$ which equals a function in $\Cinf(\R)$ in the sense of
generalized distributions need not belong to $\Ginf(\R)$ (compare
with \cite[Proposition 3.17]{NPS:98}).
\end{rem}
\begin{rem} \label{lowerordertermsmatter}
The operator $P(\xi) = a^2\xi^2 +1$, with $a$ nonnegative and invertible, is
SH-elliptic with radius $1$, since one easily verifies the following inequalities, valid when
$\eps\in(0,1/4)$:
\[
    |P_\eps(\xi)| \geq \eps^{N+1} (1 + |\xi|)^2,
    \frac{|\d_\xi P_\eps(\xi)|}{|P_\eps(\xi)|} \leq \frac{4}{1 + |\xi|},
    \frac{|\d^2_\xi P_\eps(\xi)|}{|P_\eps(\xi)|} \leq \frac{8}{(1 + |\xi|)^2}.
\]
We will see in the next section that SH-ellipticity implies $\Ginf$-hypoellipticity.
This shows that lower order terms do matter.
\end{rem}
\begin{ex} \label{ex4}
The multiplication operator $P(x)$ on $\R$ given by $P_\eps(x) = 1 + x^2/\eps$:
The zero order operator $P$ is obviously S-elliptic. It is
not WH-elliptic, because an estimate of the form
\[
   \big|\d_x^2 P_\eps(x)\big| = \frac{2}{\eps} \leq C|P_\eps(x)| = 1 + x^2/\eps
\]
as $\eps \to 0$ does not hold at $x = 0$. Further, $P$ is not $\Ginf$-hypoelliptic; the
solution $u$ to the equation $Pu = 1$ does not belong to $\Ginf(\R)$. In fact, it is given by
\[
  u_\eps(x) = \big(1 + \frac{x^2}{\eps}\big)^{-1}
    = \sum_{k=0}^\infty \big(-\frac{1}{\eps}\big)^k x^{2k}
\]
where the series representation holds for $|x| < \sqrt{\eps}$. At $x = 0$ the
derivatives are
\[
   \d_x^{2\alpha}u_\eps(0) = \frac{1}{2\alpha !}\big(-\frac{1}{\eps}\big)^\alpha
\]
and hence do not have a uniform, finite order independently of $\alpha$.
\end{ex}
\begin{rem}
Contrary to the non-regular solution discussed in Remark \ref{remarksinxovereps},
$u$ is not even equal to an element of $\Ginf(\R)$ in the sense of generalized distributions.
In fact, for $\psi \in \D(\R)$,
\begin{equation} \label{ex4computation}
   \int \big(1 + \frac{x^2}{\eps}\big)^{-1} \psi(x)\,dx
    = \sqrt{\eps} \int (1 + y^2)^{-1} \psi(\sqrt{\eps}x)\,dx \to 0
\end{equation}
as $\eps \to 0$, showing that $u$ is associated with zero. But $u$ is not zero
in the sense of generalized distributions, because if it were, the right hand side
of (\ref{ex4computation}) should be $O(\eps^q)$ for every $q \geq 0$.
This is not the case if $\psi(0) \neq 0$; then it actually has the precise
order $\kappa = \sqrt{\eps}$, since the integral converges to the finite limit $\pi \psi(0)$.

In addition, Example \ref{ex4} exhibits an invertible element of $\G(\R)$, namely
the class of $(P_\eps(x))_\eps$, which is a member of $\Ginf(\R)$ but whose
multiplicative inverse does not belong to $\Ginf(\R)$.
\end{rem}

%% file: sec5a.tex
\section{The elliptic regularity result for WH-elliptic operators with
constant coefficients}

We start this section with a general regularity result for the
constant coefficient case.
Consider an operator $P$ with symbol $P(\xi) = \sum_{|\ga| \leq m} a_\ga
\xi^\ga$, with coefficients $a_\ga \in \gC$, which is assumed to be
WH-elliptic. Thus it satisfies the condition
\begin{gather*}
    \exists N > 0\; \exists \eps_0 > 0\; \forall \eps\in
        (0,\eps_0)\; \exists r_\eps > 0 \text{\ \ such
        that}\\ |P_\eps(\xi)| \geq \eps^N (1 + |\xi|)^m
\end{gather*}
for all $\eps\in (0,\eps_0)$ and for all $\xi\in\R^n$ with $|\xi| \geq
r_\eps$, as well as an estimate ($\al\in\N_0^n$)
\[
    |\d^\al P_\eps(\xi)| \leq C_\al |P_\eps(\xi)| (1 + |\xi|)^{-|\al|}
\]
for $\eps$ and $\xi$ in the same range.

We begin by constructing a generalized parametrix for the operator $P$. Let
$\chi\in\Cinf(\R^n)$, $\chi(\xi) \equiv 0$ for $|\xi| \leq 1$ and $\chi(\xi)
\equiv 1$ for $|\xi| \geq 2$ and put
\[
    \chi_\eps(\xi) = \chi(\xi / r_\eps)
\]
and
\[
    Q_\eps = \F^{-1}(\frac{\chi_\eps}{P_\eps}),
        h_\eps = \F^{-1}(1 - \chi_\eps).
\]
It is clear that, for fixed $\eps\in(0,\eps_0)$, $Q_\eps\in\S'(\R^n)$ and
$h_\eps\in\S(\R^n)$. We also have
\begin{equation}\label{Q_parametrix}
    P_\eps(D) Q_\eps = \F^{-1}(\chi_\eps) = \de - h_\eps
\end{equation}
where $\de$ denotes the Dirac measure. The family $(Q_\eps)_{\eps\in(0,1]}$
of tempered distributions will serve as the generalized parametrix.

\begin{lemma} \label{chi_lemma} There is $N\geq 0$  such that for all
$\al\in\N_0^n$ there is $C_\al > 0$ such that
\begin{equation}\label{chi_est}
    |\d^\al (\frac{\chi_\eps(\xi)}{P_\eps(\xi)})| \leq
       C_\al \eps^{-N} (1 + |\xi|)^{-m-|\al|}
\end{equation}
for all $\xi\in\R^n$ and $\eps\in (0,\eps_0)$; $C_0$ can be chosen to
be $1$.
\end{lemma}
\begin{proof} If $|\al| = 0$ the assertion is obvious from the
hypothesis. We use induction over
$|\al|$.
>From the Leibniz rule we obtain, for $|\al| \geq 1$,
\[
  \d^\al (\chi_\eps / P_\eps) =
    \d^\al \chi_\eps / P_\eps - \sum_{\be < \al}
      \binom{\al}{\be} \d^\be (\chi_\eps / P_\eps)
        \d^{\al - \be} P_\eps / P_\eps.
\]
We have that $\d^\al \chi_\eps(\xi) \not= 0$ only if $|\xi| \leq 2
r_\eps$. In this range, $1 \leq (1 + 2 r_\eps)^{|\al|} (1 +
|\xi|)^{-|\al|}$, hence
\begin{multline*}
  | \d^\al\chi_\eps / P_\eps| \leq
     \linf{\d^\al\chi_\eps} \eps^{-N} (1 + 2 r_\eps)^{|\al|}
       (1 + |\xi|)^{- m - |\al|} \\
     = r_\eps^{-|\al|}\linf{\d^\al\chi} \eps^{-N} (1 + 2 r_\eps)^{|\al|}
       (1 + |\xi|)^{- m - |\al|}\\
     \leq C'_\al \eps^{-N} (1 + |\xi|)^{- m - |\al|}.
\end{multline*}
Here, $N$ is chosen as in the ellipticity condition (and is
independent of $\al$).
Furthermore, $|\d^{\al - \be} P_\eps(\xi) / P_\eps(\xi)|
\leq C''_{\al - \be} (1 + |\xi|)^{|\be| - |\al|}$ and
$|\d^\be (\chi_\eps(\xi) / P_\eps(\xi))| \leq
C_\be \eps^{-N} (1 + |\xi|)^{- m - |\be|}$ when
$|\be| < |\al|$ by assumption.
Thus the conclusion follows.
\end{proof}
\begin{lemma} \label{h_lemma} $\linf{\d^\al h_\eps} \leq 2^{|\al|}
r_\eps^{|\al|+n} \lone{1-\chi}$, for all $\eps\in (0,\eps_0)$. In
particular, $(h_\eps)_{\eps\in(0,1)} \in \EM^\infty(\R^n)$ if $r_\eps$
is of slow scale.
\end{lemma}
\begin{proof}
\begin{multline*}
  \linf{\d^\al h_\eps} \leq
    \int |\xi|^{|\al|} |1 - \chi_\eps(\xi)|\, d\xi\\
  \leq (2 r_\eps)^{|\al|} \int |1 - \chi(\frac{\xi}{r_\eps})|\, d\xi
  = 2^{|\al|} r_\eps^{|\al|+n} \lone{1 - \chi}.
\end{multline*}
The second inequality follows from the fact that $1 - \chi_\eps(\xi))
\equiv 0$ when $|\xi| \geq 2 r_\eps$.
\end{proof}
\begin{lemma} \label{Q_lemma} For every $K \Subset \Om$ and $s > n/2$
there is a constant $C > 0$ such that the Sobolev estimate
\[
  \linf{Q_\eps * \vphi} \leq C \eps^{-N} \norm{\vphi}{W^{s,\infty}}
\]
holds for all $\vphi \in \D(K)$ and all
$\eps\in(0,\eps_0)$.
\end{lemma}
\begin{proof}
\[
  \linf{Q_\eps * \vphi} \leq \linf{\chi_\eps / P_\eps}
    \lone{\FT{\vphi}}
  \leq C\eps^{-N} \norm{\vphi}{H^s} \leq C \eps^{-N}
    \norm{\vphi}{W^{s,\infty}}
\]
with $C^2 = \int (1 + |\xi|)^{-2s}\, d\xi$ by usual Sobolev space
arguments. The second inequality
uses the fact that $ \linf{\chi_\eps / P_\eps} \leq \eps^{-N}$ by
Lemma \ref{chi_lemma}.
\end{proof}
\begin{prop}\label{Q_prop}
$\big( Q_\eps \vert_{\R^n\setminus\{0\}} \big)_{\eps \in (0,1]}$
defines and element of $\EM^\infty(\R^n \setminus\{0\})$.
\end{prop}
\begin{proof}
Take $\varphi \in \D(\R^n), \varphi \equiv 1$ in a neighborhood of
zero, and put $\psi(x) = \varphi(\sigma x) - \varphi(x/\sigma)$ with
$\sigma > 0$. By taking $\sigma$ sufficiently small, every compact
subset of $\R^n\setminus\{0\}$ eventually lies in the region where
$\psi \equiv 1$. Thus it suffices to establish the $\EM^\infty$-estimates
(\ref{regefu}) for $\psi Q_\eps$. Its Fourier transform equals
\[
   {\cal F}(\psi Q_\eps)
   = \int \frac{\chi_\eps(\xi - \eta}{P_\eps(\xi - \eta)} \FT{\psi}(\eta)\,d\eta
   = \sum_{|\beta|=q}\frac{1}{\beta !}\int
       \d^\beta\Big(\frac{\chi_\eps}{P_\eps}\Big)(\theta)\,\eta^\beta\FT{\psi}(\eta)\,d\eta
\]
for every $q \geq 1$, where $\theta$ lies between $\xi$ and $\xi - \eta$.
This follows by Taylor expansion and observing that $\FT{\psi} \in \S(\R^n)$
has all its moments vanishing. By Lemma \ref{chi_lemma} and Peetre's inequality
we have for $|\beta| = q$ that
\[
  \big|\d^\beta\Big(\frac{\chi_\eps}{P_\eps}\Big)(\theta)\big|
    \leq C_q\,\eps^{-N}\big(1 + |\theta|\big)^{-m-q}
     \leq C_q'\,\eps^{-N}\big(1 + |\xi|\big)^{-m-q} \big(1 + |\eta|\big)^{m+q}\,.
\]
Let $\alpha \in \N_0^n$ and choose $q$ large enough so that
$|\alpha| - m - q < -n$. Then
\[
   \xi^\alpha{\cal F}\big(\psi Q_\eps\big)(\xi) \in L^1(\R^n)\,, \quad
     \d^\alpha\big(\psi Q_\eps\big)(x) \in L^\infty(\R^n)
\]
and
\[
   \linf{\d^\alpha(\psi Q_\eps)}
      \leq C_q'\,\eps^{-N} \int \big(1 + |\eta|\big)^{m+q}|\eta|^q|\FT{\psi}(\eta)|\,d\eta
       \int |\xi|^{|\alpha| - m - q}\,d\xi\,.
\]
This proves that $\psi Q_\eps \in \Cinf(\R^n)$ with all derivatives satisfying
a bound of order $\eps^{-N}$.
\end{proof}

\begin{thm}\label{ell_reg_thm}
Let the operator  $P(D) = \sum_{|\ga|\leq m} a_\ga D^\ga$,
with $a_\ga\in\gC$, be WH-elliptic with radius $r_\eps$ of slow
scale. Let $f\in\Ginf(\Om)$ and $u\in\G(\Om)$ be a solution to $P(D) u
= f$. Then $u\in\Ginf(\Om)$.
\end{thm}
\begin{proof} Let $\om \Subset \Om$ and choose
$\vphi\in\D(\Om)$, $\vphi \equiv 1$ on $\om$. Then
\[
  P(D)(\vphi u) = \vphi f + v
\]
where $v \equiv 0$ on $\om$ and $v$ has compact support in $\Om$. It
suffices to show that $\vphi u$ enjoys the $\Ginf$-property on every
compact set $K \subset \om$.

We have that (see (\ref{Q_parametrix}))
\begin{multline*}
  \vphi u_\eps = P_\eps(D) Q_\eps * (\vphi u_\eps) +
        h_\eps * (\vphi u_\eps) \\
  = Q_\eps * (\vphi f_\eps + v_\eps) + h_\eps * (\vphi u_\eps)\\
  = Q_\eps * (\vphi f_\eps) + (\psi Q_\eps) * v_\eps +
     \big((1 - \psi)Q_\eps\big) * v_\eps + h_\eps * (\vphi u_\eps)\\
  = \text{I} + \text{II} + \text{III} + \text{IV}
\end{multline*}
where $\psi$ is a cut-off supported in a small neighborhood of zero,
$\psi(0)=1$. We shall prove the $\Ginf$-property of each term on $K$.

For term (I) we have
\[
  \linf{\d^\al\big( Q_\eps * (\vphi f_\eps)  \big)} \leq
     C \eps^{-N} \norm{\d^\al(\vphi f_\eps)}{W^{s,\infty}}
\]
by Lemma \ref{Q_lemma}, when $s > n/2$. But $f\in\G(\Om)$, so the
latter term has a bound of order $\eps^{-N'}$ independently of
$\al$. Concerning term (II), we may choose the support of
the cut-off $\psi$ so small that this term actually vanishes on $K$.

Coming to term (III), we have by Proposition \ref{Q_prop} that
$(1-\psi) Q$ belongs to $\Ginf(\R^n)$. To estimate
$\d^\al\big( (1-\psi)Q_\eps \big) * v_\eps$, we let all the
derivatives fall on the first factor, observe that $v_\eps \in
\Con(\ovl{\om})$ and evoke the continuity of the convolution map
$\Con(\R^n)\times \Con(\ovl{\om}) \to \Con(\R^n)$ to conclude that
$\d^\al\big( (1-\psi)Q_\eps \big) * v_\eps$ has a bound of order
$\eps^{-N}$ independently of $\al$, uniformly on each compact subset
of $\R^n$. Finally, term (IV) is treated by the same argument,
observing that $h\in\Ginf(\R^n)$ by Lemma \ref{h_lemma}.
\end{proof}

\begin{rem} The operator $P(\xi) = - a \xi^2 +1$
from Example \ref{ex3} was shown to be WH-elliptic with radius
$r_\eps = s/a_\eps, s > 1$. From Theorem \ref{ell_reg_thm} and the
explicit solution of the homogeneous equation we may now assert
that it is $\Ginf$-hypoelliptic if and only if $r_\eps$ is
slow scale. This once again emphasizes the importance of
the slow scale property.
\end{rem}
\begin{rem}
If $P(D)$ is W-elliptic with radius $r_\eps$, then all
real roots of $P_\eps(\xi) = 0$ lie within the radius $r_\eps$. Let
$m_\eps = \max \{ |\xi| : P_\eps(\xi) = 0, \xi\in\R^n  \}$. Since
$m_\eps \leq r_\eps$, a necessary requirement for the conditions of
Theorem \ref{ell_reg_thm} to hold is that $m_\eps$ is slow
scale. If $m_\eps$ is moderate, the slow scale property is also
necessary for the solutions of $P(D) u = 0$ to belong to $\Ginf$, as
is seen from the solutions $u_\eps(x) = \exp(i x \xi_\eps)$ where
$P_\eps(\xi_\eps) = 0$.
\end{rem}
The remainder of this section is devoted to second order operators
\[
    P(D) = \sum_{i,j=1}^n a_{ij} D_i D_j + \sum_{j=1}^n b_j D_j + c\,,
\]
in which case we have a regularity result under somewhat different
assumptions than in Theorem \ref{ell_reg_thm}.
For its generalized constant coefficients, we assume that $a_{ij} \in \gR$,
$b_j, c \in \gC$. Such an operator is called \emph{G-elliptic} (for
\emph{generalized elliptic}), if the matrix $A = (a_{ij})_{i,j} \in
\gR^{n^2}$ is symmetric and positive definite. Thus all eigenvalues
$\la_1$, $\ldots$, $\la_n \in \gR$ of $A$ are invertible and nonnegative.
Employing $\R$-linear algebra on representatives at fixed $\eps > 0$, one sees
that there is an orthogonal matrix $Q$ with coefficients in $\gR$ such that
\[
    A = Q^T \Lambda Q, \quad \Lambda = \text{diag}(\la_1,\ldots,\la_n).
\]
Define
\[
    \be_i = \sum_{j=1}^n Q_{ij} b_j,
\]
and, putting $\ga = c - \sum_{j=1}^n \be_j^2 / 4 \la_j$,  let
\[
    \tilde{\la}_i = \begin{cases}
                    \la_i & \text{if $\ga = 0$},\\
                    \la_i / |\ga| & \text{if $\ga$ is invertible}.
                    \end{cases}
\]
The case of non-zero, non-invertible $\gamma$ is outside the scope of the
following result.

\begin{prop}\label{2nd_order_prop} Let $P(D)$ be a second order
G-elliptic operator and assume that $\ga$ is either equal to zero or
invertible. Let $f\in\Ginf(\R^n)$ and let $u\in\G(\R^n)$ be a solution to
$P(D) u = f$. If $\be_j / \la_j$ is of log-type and $\tilde{\la}_i$,
$\tilde{\la}_i^{-1}$ is slow scale, $i,j=1,\ldots,n$, then $u\in\Ginf(\R^n)$.
\end{prop}
\begin{proof} The proof proceeds by stepwise reduction to the Laplace- or
Helmholtz equation. First, we change the independent variable to $y = Q x$ and
define $v(y) = u(Q^T y)$, $g(y) = f(Q^T y)$, that is, we define
representatives $v_\eps(y) = u_\eps(Q_\eps^T y)$ and $g_\eps(y) = f\eps(Q_\eps^T y)$.
The columns of $Q$ and
$Q^T$ have length $1$, so $v$ is moderate iff $u$ is moderate, and
$u\in\Ginf(\R^n)$ iff $v\in\Ginf(\R^n)$; the same holds of $f$ and $g$.
Further,
\[
    P(D) u = f
\]
is equivalent with
\begin{equation}\label{erstred}
    \sum_{i=1}^n \la_i D_i^2 v + \sum_{j=1}^n \be_j D_j v + c v = g
\end{equation}
in $\G(\R^n)$ with respect to the new variables $y =Q x$. We may rewrite
(\ref{erstred}) to
\begin{equation}\label{zweitred}
    \sum_{i=1}^n \la_i \big(D_i + \frac{\be_i}{2\la_i}\big)^2 v +
        \big(c - \sum_{j=1}^n \frac{\be_j^2}{4\la_j}\big) v = g
\end{equation}
and put
\[
    w(y) = \prod_{j=1}^n \exp\big(-i\frac{\be_j}{2\la_j} y_j\big) v(y).
\]
If $\be_j/\la_j$ is of log-type, this transformation respects moderateness,
and so (\ref{zweitred}) is equivalent with
\[
    \sum_{i=1}^n \la_i D_i^2 w(y) + \ga w(y) = h(y)
\]
in $\G(\R^n)$, where $h(y) = \prod_{j=1}^n\exp(-i\frac{\be_j}{2\la_j} y_j) g(y)$.
Further, if $\be_j/\la_j$
is of log-type, then $v\in\Ginf(\R^n)$ iff $w\in\Ginf(\R^n)$, and
similarly for $g$ and $h$.

Let  $\tilde{h}(y) = h(y) / |\ga|$ if $\ga$ is invertible, and
$\tilde{h}(y) = h(y)$ if $\ga = 0$. According to the hypotheses of the
proposition, we may rewrite the latter equation as
\[
    \sum_{i=1}^n \tilde{\la}_i D_i^2 w(y) + \sig w(y) = \tilde{h}(y)
\]
where $\sig = 0$ or else $|\sig| = 1$. Finally, we put
\begin{eqnarray*}
    \tilde{w}(y_1,\ldots,y_n) & = &
        w(\sqrt{\tilde{\la}_1} y_1,\ldots,\sqrt{\tilde{\la}_n} y_n),\\
    \tilde{\tilde{h}}(y_1,\ldots,y_n) & = &
        \tilde{h}(\sqrt{\tilde{\la}_1} y_1,\ldots,\sqrt{\tilde{\la}_n} y_n),
\end{eqnarray*}
and arrive at the equation
\[
    -\Delta \tilde{w}(y)+ \sig \tilde{w}(y) = \tilde{\tilde{h}}(y).
\]
As above,  the hypotheses imply that $w\in\Ginf(\R^n)$ iff
$\tilde{w}\in\Ginf(\R^n)$, and the same holds for $\tilde{h}$ and
$\tilde{\tilde{h}}$.

Collecting everything, we see that $f\in\Ginf(\R^n)$ iff
$\tilde{\tilde{h}}\in\Ginf(\R^n)$. But the operator $|\xi|^2 + \sig$ is clearly
SH-elliptic (with radius r = 1), so Theorem \ref{ell_reg_thm}
shows that $\tilde{w}\in\Ginf(\R^n)$ which in turn
implies that $u\in\Ginf(\R^n)$ as desired.
\end{proof}

\begin{ex}\label{2orderexample} The second order homogeneous ODE
\[
    \big(\la \frac{d^2}{dx^2} + b \diff{x} + c\big) u(x) = 0
\]
has the solution
\[
   u(x) = C_1 \exp (\mu_+x) + C_2 \exp (\mu_- x)
\]
where
\[
\mu_{\pm} =
-\frac{b}{2\la} \pm \frac{1}{\sqrt{\la}}\sqrt{\frac{b^2}{4\la} - c}\,.
\]
$(u_\eps)_\eps$ is moderate and thus defines an element of $\G(\R)$ if either $\mu_{\pm}$
has nonzero real part and is of log-type or else if its real part is zero and it
is slow scale. In both cases, the solution belongs to $\Ginf(\R)$. This
illustrates the hypotheses required in Proposition \ref{2nd_order_prop}.
\end{ex}

\begin{rem} The conditions stated in Theorem
\ref{ell_reg_thm} and Proposition \ref{2nd_order_prop} are independent (for
second order G-elliptic operators) as can be seen from the following examples:

First, the operator $P(\xi) = a^2 \xi_1^2 + \xi_2^2$ discussed in Example \ref{ex2}, with
$a_\eps = 1 / |\log(\eps)|$, is not WH-elliptic, so Theorem \ref{ell_reg_thm} does
not apply, but Proposition \ref{2nd_order_prop} does, and so the operator is $\Ginf$-hypoelliptic.

Second, the operator $P(\xi) = a^2\xi^2 + 1$ was seen to be SH-elliptic in
Remark \ref{lowerordertermsmatter}, so it is $\Ginf$-hypoelliptic by
Theorem \ref{ell_reg_thm}. However, if $a_\eps = \eps$ or any positive power
thereof, it does not satisfy the log-type property required in
Proposition \ref{2nd_order_prop}. The corresponding homogeneous
differential equation $(-\eps^2 \d_x^2 +  1) u_\eps(x) = 0$
has the solution $u_\eps(x) = C_1 \exp(x/\eps) + C_2 \exp(-x/\eps)$, which
does not yield an element of $\Ginf(\R)$. This does not contradict the
$\Ginf$-hypoellipticity of the operator $P(D)$, because $(u_\eps)_\eps$ is not
moderate and so does not represent a solution in $\G(\R)$.
We will take up this observation in Example \ref{nonsolvableop} below.

If the coefficient $\gamma$ entering the hypotheses of Proposition \ref{2nd_order_prop}
is a zero divisor, the operator $P(\xi) = \xi^2 + \gamma$
may or may not be $\Ginf$-hypoelliptic. It is so, if $\gamma \geq 0$ (similar to
Remark \ref{lowerordertermsmatter}). It is not, if $\gamma_\eps \in \{0, -1/\eps^2\}$ in
an interlaced way as $\eps \to 0$ (following Example \ref{ex3}).
\end{rem}
\begin{ex} \label{nonsolvableop}
The regularity result of Theorem \ref{ell_reg_thm}
can be used to prove nonexistence of solutions: For example, let
$a \in \gR$ be the class of $a_\eps = \eps$. Then the homogeneous equation
\begin{equation} \label{nonsolvableeq}
   \big(-a^2 \frac{d^2}{dx^2} +  1\big) u(x) = 0
\end{equation}
has no nontrivial solution in $\G(\Omega)$ on whatever open subset $\Omega \subset \R$. To see this,
assume there is a solution, with representative $(u_\eps)_\eps$, say. Then
\[
   \big(-\eps^2 \frac{d^2}{dx^2} +  1\big) u_\eps(x) = n_\eps(x)
\]
for some $(n_\eps)_\eps \in \NN(\Omega)$. It follows that
\[
   u_\eps'' = \frac{1}{\eps^2}\big(u_\eps - n_\eps\big)\,,\quad
   u_\eps^{(4)} = \frac{1}{\eps^2}\big(u_\eps'' - n''_\eps\big)
       = \frac{1}{\eps^4}\big(u_\eps - n_\eps\big) - \frac{1}{\eps^2} n''_\eps\,,
\]
and so on, so that
\[
    \big(u_\eps^{(2\alpha)} - \frac{1}{\eps^{2\alpha}} u_\eps\big)_\eps \in \NN(\Omega)
\]
for every $\alpha \in \N$. On the other hand, the operator in (\ref{nonsolvableeq})
is SH-elliptic, hence $\Ginf$-hypoelliptic. Therefore, there is $N \geq 0$ such that
$u_\eps^{(2\alpha)} = O(\eps^{-N})$ for every $\alpha \in \N$. It follows that
$u_\eps = O(\eps^{-N+2\alpha})$ for every $\alpha \in \N$, whence $(u_\eps)_\eps \in \NN(\Omega)$
and $u = 0$ in $\G(\Omega)$. 
Observe that the nonexistence result depends crucially on the asymptotic
properties of the generalized coefficient $a^2$. If we take $a\in\gR$
invertibel such that $1/a$ is of logarithmic type then equation
\ref{nonsolvableeq} is nontrivially solvable (see Example \ref{2orderexample}).
\end{ex}





%% file: sec6a.tex

\section{Micro
-elliptic regularity for first order operators with variable coefficients of
 slow scale}

For partial differential operators with smooth coefficients the elliptic
regularity can be considered a special instance of microlocal
non-characteristic regularity. This is expressed in terms of a general
relation combining the wave front sets of the distributional solution and
of the right-hand side in the PDE with the characteristic set of the
operator. It states that for any partial differential operator $P(x,D)$
with $\Cinf$ coefficients and distribution $u$ we have the following
inclusion relation:
\begin{equation}\label{WF_rel}
     \WF(u) \subseteq \WF(Pu) \cup \Char(P)
\end{equation}
(for subsets of $\CO$, the cotangent space over $\Om$ with the zero section
removed.) Recall that $\Char(P) = P_m^{-1}(0) \cap
\CO$ and thus depends only on the principal symbol of the operator.

The concept of microlocal regularity of a Colombeau function follows the
classical idea of adding information about directions of rapid decrease in the
frequency domain upon localization in space (cf.\
\cite{Hoermann:99,HK:01,NPS:98}). It refines $\Ginf$ regularity in the same
way as the distributional wave front set does with $\Cinf$ regularity, i.e.,
the projection of the (generalized) wave front set into the base space equals
the (generalized) singular support.

We briefly recall the definition of the generalized wave front set of a
Colombeau function: $u\in\G(\Om)$ is said to be \emph{microlocally regular} at
$(x_0,\xi_0)\in\CO$ if (for a representative $(u_\eps)_\eps\in\EM(\Om)$) there
is an open neighborhood $U$ of $x_0$ and a conic neighborhood $\Ga$ of $\xi_0$
such that for all $\vphi\in\D(U)$ it holds that $\F(\vphi u)$ is rapidly
decreasing in $\Ga$, i.e., $\exists N\in\R$ $\forall l\in\N_0$ $\exists C > 0$
$\exists \eps_0 > 0$:
\begin{equation}\label{rap_dec}
  |(\vphi u_\eps)\FT{\ }(\xi)| \leq C
  \eps^{-N} (1 + |\xi|)^{-l} \qquad \forall\xi\in\Ga, \forall\eps\in (0,\eps_0).
\end{equation}
(Here, $\FT{\ }$ denotes classical Fourier transform and $\F(\vphi u)$ the
corresponding generalized Fourier transform of the compactly supported
Colombeau function $\vphi u$.) The
\emph{generalized wave front set} of $u$, denoted $\WF_g(u)$, is defined as
the complement (in $\CO$) of the set of pairs $(x_0,\xi_0)$ where $u$ is
microlocally regular.
\begin{rem} \label{rem_5_1} 
(i) Since $\F(\vphi u)$ is temperate it suffices to require an estimate of
 the form (\ref{rap_dec}) for all $\xi\in\Ga$ with $|\xi| \geq r_\eps$ where
 $(r_\eps)_\eps$ is of slow scale. Indeed, there are $M\in\R$ and $\eps_1 > 0$
 such that \[ |(\vphi u_\eps)\FT{\ }(\xi)| \leq \eps^{-M} (1 + |\xi|)^{M}
 \qquad \forall\xi\in\R^n, \forall\eps\in (0,\eps_1). \] When $|\xi| \leq
 r_\eps$ the right-hand side is bounded as follows: \[ \eps^{-M} (1 +
 |\xi|)^{M} \leq \eps^{-M} (1 + r_\eps)^{M+l} (1 + |\xi|)^{-l} \leq
 \eps^{-M-1} (1 + |\xi|)^{-l}. \] Hence taking $\max(M+1,N)$,
 $\min(\eps_0,\eps_1)$ as new $N$, $\eps_0$ one arrives at (\ref{rap_dec})
 valid for all $\xi\in\Ga$.

(ii) If $v\in\G_c(\Om)$ and $(r_\eps)_\eps$ is not of slow scale then
 the  rapid decrease property (\ref{rap_dec}) does not follow from the
 corresponding estimates in the regions $|\xi| \geq r_\eps$. Consider, for
 example, a model delta net $v_\eps = \rho(./\eps)/\eps$ where
 $\rho\in\D(\R)$. Then $\FT{v_\eps}$ is not rapidly decreasing in $\xi > 0$
 and $\xi < 0$ but
\begin{multline*}
    |\FT{v_\eps}(\xi)| =
    |\FT{\rho}(\eps \xi)| = |\int e^{-i\eps x \xi} \rho(x)\, dx| =
        (\eps|\xi|)^{-l} |\int e^{-i\eps x \xi} \rho^{(l)}(x)\, dx| \\
    \leq C_l (\eps|\xi|)^{-l} = C_l (\eps\sqrt{|\xi|})^{-l} |\xi|^{-l/2}
    \leq C_l |\xi|^{-l/2}
    \qquad \text{if } |\xi| \geq \frac{1}{\eps^2} =: r_\eps.
\end{multline*}
\end{rem}

In view of the examples in Section 4 one cannot expect to obtain an
extension of (\ref{WF_rel}) to arbitrary operators with $\Ginf$-coefficients
by designing a notion of generalized characteristic set based solely
on the principal part.
In general, the lower order terms of the symbol do have an effect on the
regularity properties, even for smooth principal part: consider
the symbol $P_\eps(\xi) = \xi - 1/\eps$ whose corresponding operator admits
the non-regular solution $u_\eps(x) = \exp(i x/\eps)$ to the homogeneous
equation. However, in case of first order operators with variable coefficients
a direct approach shows sufficiency of two further assumptions to restore a
microlocal regularity relation of the type (\ref{WF_rel}). Both are
requirements of slow scale: One about regularity of the coefficients and the
other in terms of lower bounds on the principal symbol over conic regions.

An identification of adequate conditions in the general case, in particular,
finding an appropriate notion of characteristic set of the operator, remains
open. As suggested by the results and examples of Sections 3 and 4 the latter
might have to include the influence of lower order terms in an essential way.

We first introduce an auxiliary notion to replace $\Char(P)$ in our variant of relation
(\ref{WF_rel}) for first order operators.
\begin{defn}\label{Wsc_def}
Let $P(x,D)$ be a partial differential operator of order $m$ with coefficients
in $\G(\Om)$ and let $(x_0,\xi_0)\in\CO$. $P$ is said to be \emph{W-elliptic
with slow scales} at $(x_0,\xi_0)$,
\emph{\Wsc-elliptic} in short, if for some representative
$(P_\eps(x,\xi))_\eps$ the following holds: there is an open
neighborhood $U$ of $x_0$, a conic neighborhood $\Ga$ of $\xi_0$,
slow scale nets $(s_\eps)_\eps$, $(r_\eps)_\eps$ with $s_\eps > 0$, $r_\eps >
0$, and $\eps_0 > 0$ such that
 \begin{equation}\label{Wsc_ell}
    |P_\eps(x,\xi)| \geq s_\eps^{-1} (1 + |\xi|)^m
 \end{equation}
for all $x\in U$, $\xi\in\Ga$, $|\xi| \geq r_\eps$, $\eps \in (0,\eps_0)$.
We denote by $\WscEll(P)$ 
the subset of pairs in $\CO$ where $P$ is \Wsc-elliptic.
\end{defn}

\begin{rem}\label{not_rem}
(i) In the theory of propagation of singularities for operators with
smooth coefficients it is the
complement of the ellipticity set, namely the characteristic set $\Char(P)$,
which plays a dominant role. In view of the remarks at the beginning of this
section and the variety of ellipticity notions used in the Colombeau context
so far we refrain from introducing (yet another) notion of
\emph{generalized characteristic set}. Finding an appropriate definition for a
theory of propagation of singularities for operators with nonsmooth
coefficients is the subject of ongoing and future research. For the purpose of
the present paper we prefer the notation $\WscChar{P}$, the
complement $\WscChar{P}$ of $\WscEll(P)$ in $\CO$.
(ii) Clearly, \Wsc-ellipticity implies $W$-ellipticity. The symbol
$P_\eps(\xi) = \eps \xi$, which is $SH$-elliptic, shows that the converse does
not hold.
\end{rem}

\begin{ex}\label{Char_ex} 
(i) $P_\eps(\xi) = \xi - 1/\eps$ gives $\WscChar{P} \- = \R\times\R\zs$
 but  $\WscChar{P_1} = \emptyset$, whereas $Q_\eps(\xi) = \xi - \log(1/\eps)$
 yields
 $\WscChar{Q} = \emptyset$. Slightly more general, let $(r_\eps)_\eps$,
 $(s_\eps)_\eps$ be slow scale nets, $|r_\eps| > 0$, then
 $P_\eps(\xi) = r_\eps^{-1} \xi - s_\eps$ is \Wsc-elliptic at every
 $(x_0,\xi_0)\in \R\times\R\zs$.

(ii) For the operator $P_\eps(\xi) = \eps \xi_1 + i \xi_2$ (from
 Example \ref{ex1})
 we obtain $\WscChar{P} = \R^2 \times (\R\times \{0\})\zs$
 since (\ref{Wsc_ell}) is valid with constant $s_\eps$ and $r_\eps$ in any
 cone $|\xi_2| \geq c |\xi_1| > 0$. We observe that the wave front set of the
 solution $u$ to $Pu = 0$ is a subset of $\WscChar{P}$. Indeed, we have
 $u_\eps(x_1,x_2) = \exp(i x_1 /\eps)\cdot \exp(- x_2)$ and application of
 \cite[Lemma 5.1]{HK:01} proves the inclusion.

 Note that $\singsupp_g(u) = \R^2$ (by direct inspection of derivatives) but
 the inclusion $\WF_g(u) \subset \WscChar{P}$ is nevertheless strict since
 $f_\eps(x_1) = \exp(i x_1 /\eps)$ has the `half-sided' wave front set
 $\WF_g(f) = \R\times\R_+$. (To see this one easily checks that
 $(\vphi \cdot \exp(i./\eps))\FT{\ }(\xi) = \FT{\vphi}(\xi - 1/\eps)$ is
 rapidly decreasing if and only if $\xi < 0$.)

(iii) Let $P_\eps(x,y,\d_x,\d_y) = \d_y - a_\eps(x,y) \d_x$ where
 $a_\eps\in\EM(\R^2)$ is real-valued and bounded uniformly with respect to
 $\eps$. Put $c_1 = \inf_{x,y,\eps} a_\eps$, $c_2 = \sup_{x,y,\eps} a_\eps$
 then any pair $((x_0,y_0),(\xi_0,\eta_0)) \in \R^2\times\R^2\zs$ with $\eta_0
 \not\in [c_1,c_2] \cdot \{\xi_0\}$ is a point of \Wsc-ellipticity of $P$.
\end{ex}

The second slow scale condition used in the theorem to follow is introduced as
a strong regularity property of the prospective coefficients in the operators.
\begin{defn}\label{sc_fun}
$v\in\G(\Om)$ is said to be of \emph{slow scale} if it has a representative
$(v_\eps)_\eps\in\EM(\Om)$ with the following property: $\forall K
\Subset \Om$ $\forall \al\in\N_0^n$ $\exists$ slow scale net
$(r_\eps)_\eps$, $r_\eps > 0$ $\exists \eps_0 > 0$ such that
 \begin{equation}
    |\d^\al v_\eps(x)| \leq r_\eps \qquad \forall x\in K, \forall \eps\in
        (0,\eps_0).
 \end{equation}
\end{defn}
\begin{rem} 
(i) Any Colombeau function of slow scale is in $\Ginf$ but clearly the
    converse does not hold.

(ii) Examples of functions of slow scale are obtained by logarithmically
    scaled embeddings: If $v\in\S'$, $\rho\in\S$ is a mollifier and we put
    $\rho^\eps(x) := (\log(1/\eps))^n \cdot \rho(\log(1/\eps)x)$  then $v_\eps
    := v * \rho^\eps$ defines the Colombeau function $v = [(v_\eps)_\eps]$
    which is of slow scale. This occurs in applications, e.g., when considering
    hyperbolic PDEs with discontinuous or nonsmooth coefficients
    (\cite{HdH:01,LO:91,O:89}).
\end{rem}

\begin{thm} \label{WF_thm} Let $P(x,D)$ be a first order partial differential
    operator with coefficients of slow scale. Then we have for any $u\in\G(\Om)$
 \begin{equation} \label{WFg_rel}
    \WF_g(u) \subseteq \WF_g(P u) \cup \WscChar{P_1}.
 \end{equation}
\end{thm}
\begin{rem} 
(i) We point out that the zero order terms of the symbol do not appear in the
 determination of $\WscChar{P_1}$. In this respect (\ref{WFg_rel}) is closer
 to the classical relation (\ref{WF_rel}) than can be expected in more
 general situations. If $\WscChar{P_1} = \emptyset$ this can be considered a
 special version of an elliptic regularity result not covered by Theorem
 \ref{ell_reg_thm} above. 

(ii) In Example \ref{Char_ex}, (ii) above the inclusion relation
 (\ref{WFg_rel}) is strict. On the other hand, as trivial examples like $P =
 1$ and $P = \diff{x}$ show we may also have equality in (\ref{WFg_rel}).
\end{rem}

The proof of Theorem \ref{WF_thm} will be based on an integration by parts
technique as in the classical regularization of oscillatory integrals. In case
of Colombeau functions regularity is coupled to the asymptotic behavior with
respect to $\eps$. Thus we have to carefully observe the interplay of this
parameter with the spatial variables when estimating Fourier integrals. The
following auxiliary result will be useful in this context.

\begin{lemma}\label{uni_lem}
Let $v\in\G(\Om)$ be microlocally regular at $(x_0,\xi_0)\in\CO$ with conic
neighborhood $U\times\Ga$ such that (\ref{rap_dec}) holds. Let $M$ be a set
and assume that $(g^\nu)_{\nu\in M}$ is a family of Colombeau functions
$g^\nu\in\G(\Om)$ with the following properties:
\begin{enumerate}
\item $\exists S \Subset U$ $\forall \nu\in M$:
    $\supp_g(g^\nu) \subseteq S$.
\item $(g^\nu)_{\nu\in M}$ satisfies a uniform $\Ginf$-property, i.e.,
    $\forall K  \Subset \Om$ $\forall \al\in\N_0^n$ $\exists q \geq 0$
$\exists \eps_0 > 0$ such that
    \[
        |\d^\al g^\nu_\eps(x)| \leq \eps^{-q} \qquad \forall x\in K, \forall
        \nu\in M, \forall \eps\in (0,\eps_0).
    \]
\end{enumerate}
Then $\F(g^\nu v)$ is rapidly decreasing in $\Ga$ uniformly with respect to
$\nu\in M$. To be more precise, for any choice of $\psi\in\D(U)$ with $\psi =
1$ on $S$ there is a real number $N'$, depending only on $\psi$ and $v$, such
that $\F(g^\nu v)$ satisfies (\ref{rap_dec}) with $C=1$ and uniform
$\eps$-power $-N'$ for all $\nu\in M$.
\end{lemma}

\begin{rem} Note that condition (ii) is satisfied in particular if each
derivative of $\g^\nu$ has slow scale bounds uniformly with respect to
$\nu\in M$.
\end{rem}

\begin{proof}
The idea of the following proof is to view this as a special case of
Theorem 3.1 in \cite{HK:01} when $\nu$ is fixed. However, to determine the
precise $\eps$-growth we have to refine the estimates along the way
appropriately.

First, we note that if $\Ga_0$ is a closed conic neighborhood of $\xi_0$ such
that $\Ga_0 \subseteq \Ga \cup \{ 0 \}$ then one can find $c > 0$ with the
following property: $\xi\in\Ga_0$, $\eta\in\Ga^c$ $\Rightarrow$ $|\xi - \eta|
\geq c |\eta|$. (See \cite[proof of Lemma 8.1.1]{Hoermander:V1} or \cite[Lemma
3.1 (i)]{HK:01} for details.)

Let $\xi\in\Ga_0$ and consider
\[
    |(g^\nu_\eps v_\eps)\FT{\ }(\xi)| = |(g^\nu_\eps \psi v_\eps)\FT{\ }(\xi)|
        = |\FT{g^\nu_\eps} * (\psi v_\eps)\FT{\ } (\xi)| \leq
            \int |\FT{g^\nu_\eps}(\xi - \eta)| |(\psi v_\eps)\FT{\ }(\eta)| \,
                d\eta.
\]
We split the integration into two parts according to the cases $\eta\in\Ga$
and $\eta\in\Ga^c$.

From the exchange formula $\eta^\al \FT{g^\nu_\eps}(\eta) = \FT{(D^\al
g^\nu_\eps)}(\eta)$ we deduce that the (global) rapid decrease estimates for
$\FT{g^\nu_\eps}$ are uniform with respect to $\nu$. This and condition (ii)
yields that $\forall p\in\N_0$ $\exists \eps_0 > 0$:
\begin{equation}\label{5_6}
 |\FT{g^\nu_\eps}(\zeta)| \leq \eps^{-q} (1 + |\zeta|)^{-p} \qquad
  \forall\zeta\in\R^n, \forall \eps\in (0,\eps_0), \forall \nu\in M.
\end{equation}

In integrating over $\eta\in\Ga$ we have rapid decrease of $|(\psi
v_\eps)\FT{\ }(\eta)|$ and hence the integrand is bounded as follows for some
$N\in\N_0$ and for all $l, k \in \N_0$ and some $\eps_0 > 0$:
\[
    |\FT{g^\nu_\eps}(\xi - \eta)| |(\psi v_\eps)\FT{\ }(\eta)| \leq
        \eps^{- q - N} (1 + |\xi - \eta|^2)^{-k/2} (1 + |\eta|^2)^{-l/2}
        \qquad \forall \eta\in\Ga, \forall\eps\in(0,\eps_0).
\]
Peetre's inequality gives $(1 + |\xi - \eta|^2)^{-k/2} \leq 2^{k/2} (1 +
|\xi|^2)^{-k/2} (1 + |\eta|^2)^{k/2}$ and we obtain for all $k\in\N_0$
\begin{multline}\label{5_7}
  \int_\Ga |\FT{g^\nu_\eps}(\xi - \eta)|
        |(\psi v_\eps)\FT{\ }(\eta)| \,d\eta \\
 \leq \eps^{-N-q}\, 2^{k/2} (1 + |\xi|^2)^{-k/2}
   \int (1 + |\eta|^2)^{(k-l)/2}\,d\eta \leq \eps^{-N-q-1} (1 + |\xi|)^{-k}
\end{multline}
when $\eps\in (0,\eps_1)$, $\eps_1$ small enough, and for all
$\nu\in M$, if $l > k + n$.

In the integral over $\eta\in\Ga^c$ we use the facts that
$(\psi v_\eps)\FT{\ }$ is temperate (in the Colombeau sense) and the cone
inequality $|\xi - \eta| \geq c |\eta|$. There is $M\in\N_0$ and $\eps_2 > 0$
such that
\begin{equation}\label{5_8}
    |(\psi v_\eps)\FT{\ }(\eta)| \leq \eps^{-M} (1 + |\eta|^2)^{M/2} \qquad
    \forall\eta\in\R^n, \forall \eps\in (0,\eps_2),
\end{equation}
which in combination with (\ref{5_6}) gives the following upper bound of the
integrand for $k, l \in\N_0$ arbitrary, $\eps \in (0,\eps_3)$, $\eps_3$ small
enough, and some $C' > 0$:
\begin{multline*}
    |\FT{g^\nu_\eps}(\xi - \eta)| |(\psi v_\eps)\FT{\ }(\eta)| \leq
        \eps^{-q - M} (1 + |\xi - \eta|^2)^{-k/2 - l/2} (1 + |\eta|^2)^{M/2} \\
        \leq \eps^{-q - M} C' (1 + |\xi|)^{-k} (1 + |\eta|)^{M + k - l}
\end{multline*}
where we have applied Peetre's inequality to the factor  $(1 + |\xi -
 \eta|^2)^{-k/2}$ and the cone inequality in estimating the factor $(1 + |\xi -
 \eta|^2)^{-l/2}$. Thus we obtain
\begin{multline}\label{5_9}
    \int_{\Ga^c} |\FT{g^\nu_\eps}(\xi - \eta)|
        |(\psi v_\eps)\FT{\ }(\eta)| \,d\eta \\
    \leq C' \eps^{-M-q} (1 + |\xi|)^{-k} \int (1 + |\eta|)^{M + k - l} \, d\eta
        \leq \eps^{-M-q-1} (1 + |\xi|)^{-k}
 \end{multline}
if $l > M + k + n$, $\eps\in (0,\eps_4)$, $\eps_4$ small enough, and for all
 $\nu\in M$. Combining (\ref{5_7}) and (\ref{5_9}) we have shown that for all
 $k\in\N_0$, $\eps_5 := \min(\eps_1,\eps_4)$
\[
    |(g^\nu_\eps v_\eps)\FT{\ }(\xi)| \leq
        (1 + |\xi|)^{-k} (\eps^{-N-q-1} + \eps^{-M-q-1})
        \qquad \forall \xi\in\Ga_0, \forall \nu\in M,
        \forall \eps\in (0,\eps_5).
\]
Since $\Ga_0$ was an arbitrary closed conic neighborhood of $\xi_0$ in $\Ga
\cup \{ 0\}$ we may put $N' = \max(N,M) + q + 1$ and the Lemma is proved.
\end{proof}

\emph{Proof of Theorem \ref{WF_thm}.} Let $(x_0,\xi_0)\in\CO$ be in the complement of the
right-hand side of (\ref{WFg_rel}) and choose $U\ni x_0$ open, $\Ga\ni \xi_0$
a conic and closed neighborhood such that both defining properties of
$\WF_g(Pu)^c$ as well as of $\WscEll(P_1)$ are satisfied when $(x,\xi)\in
U\times\Ga$. Let $\vphi\in\D(U)$ and denote by $\xib$ the projection of $\xi
\not= 0 $ onto the sphere $S^{n-1}$. We will show that the function $t \mapsto
(\vphi u_\eps)\FT{\ }(t \xib)$ is rapidly decreasing when $t \geq
\max(1,r_\eps)$ uniformly with respect to $\xib\in\Ga \cap S^{n-1}$. This
suffices to prove the theorem by Remark \ref{rem_5_1}(i).

Let $P_\eps(x,D) = \sum_{j=1}^n a_j^\eps(x) D_j + a_0^\eps(x)$ and observe that $P_{1\eps}(x,D)
\exp(-ix\xi) = - \exp(-ix\xi) P_{1,\eps}(x,\xi)$. This suggests to define the first order
differential operator $L_{\eps,\xi}(x,D)$ with parameters $\eps$, $\xi$ by
 \begin{multline}\label{L_op}
    L_{\eps,\xi}(x,D) := \Big( \frac{-1}{P_{1,\eps}(x,\xi)} P_{1,\eps}(x,D) \Big)^t \\
        = \frac{1}{P_{1,\eps}(x,\xi)} P_{1,\eps}(x,D) +
            \sum_{j=1}^n D_j(\frac{a_j^\eps(x)}{P_{1,\eps}(x,\xi)})
        = \frac{1}{P_{1,\eps}(x,\xi)} P_{1,\eps}(x,D) + q_0^\eps(x,\xi).
 \end{multline}
To avoid heavy notation in the calculations below we will henceforth denote the operator
$L_{\eps,\xi}(x,D)$ simply by $L$.

Let $\xi = t \xib$ with $t \geq \max(1,r_\eps)$. Note that $P_{1,\eps}(x,\xi)
= P_{1,\eps}(x,\xibR)\cdot t/r_\eps$, $q_0^\eps(x,\xi) = q_0^\eps(x,\xibR)
\cdot r_\eps/t$ and $1/P_{1,\eps}(x,\xibR)$ as well as $q_0^\eps(x,\xibR)$
satisfy slow scale estimates in $x$ uniformly with respect to $\xib\in \Ga\cap
S^{n-1}$ (in the sense of (ii) in Lemma \ref{uni_lem}).

Integrating by parts and applying the Leibniz rule for $P_{1,\eps}$ we have
\begin{multline*}
    (\vphi u_\eps)\FT{\ }(\xi) = \int e^{-ix\xi} L(\vphi u_\eps)(x) \, dx \\
    = \int e^{-ix\xi} \Big( \frac{1}{P_{1,\eps}(x,\xi)}
        \big( P_{1,\eps}\vphi(x) \cdot u_\eps(x) + \vphi(x)\cdot
            P_{1,\eps} u_\eps(x) \big) + q_0^\eps(x,\xi) \vphi(x) u_\eps(x)  \Big)\, dx.
\end{multline*}
We rewrite the middle term using $P_{1,\eps}u_\eps = P_\eps u_\eps - a_0^\eps
u_\eps$ and obtain
\begin{multline*}
    (\vphi u_\eps)\FT{\ }(\xi) =
        \int e^{-ix\xi} \Big(
            \frac{P_{1,\eps}\vphi(x) - a_0^\eps(x) \vphi(x)}{P_{1,\eps}(x,\xi)}
            + q_0^\eps(x,\xi) \vphi(x) \Big)
            \cdot u_\eps(x) \, dx \\
        + \frac{r_\eps}{t} \int e^{-ix\xi} \frac{\vphi(x)}{P_{1,\eps}(x,\xib)}
            P_\eps u_\eps(x) \, dx =: I_1^\eps(\xi) +
        \frac{r_\eps}{t} h_1^\eps(\xib,\xi).
\end{multline*}
The factor within parentheses in the first integral is $L\vphi(x) - a_0^\eps(x)\vphi(x) /
P_{1,\eps}(x,\xi)$ and will be abbreviated as $\vphi_{1,\eps(x,\xi)}$.

Choose $\psi\in\D(U)$ with $\psi = 1$ on $S := \supp(\vphi)$. We put
$g_\eps^{\xib}(x) := \vphi(x) / P_{1,\eps}(x,\xib)$ and observe that this
defines a family of functions satisfying properties (i) and (ii) of Lemma
\ref{uni_lem} (with $M = \Ga \cap S^{n-1}$). Since $h_1^\eps(\xib,\xi) =
(g_\eps^{\xib} P_\eps u_\eps)\FT{\ }(\xi)$ we deduce from the same lemma (with
$v = Pu$) that $\eta \mapsto h_1^\eps(\xib,\eta)$ is rapidly decreasing when
$\eta\in\Ga$, with uniform $\eps$-power, say $-K$, when $\xib$ varies in
$\Ga\cap S^{n-1}$. In particular, $t \xib \in\Ga$ and hence
$t \mapsto h_1^\eps(\xib,t \xib)$ enjoys the same decrease estimate.

In $I_1^\eps(\xi) = \int \exp(-ix\xi) \vphi_{1,\eps}(x,\xi) u_\eps(x)\, dx$ we have
$\vphi_{1,\eps}(x,\xi) = \vphi_{1,\eps}(x,\xibR)\cdot r_\eps/t$ and the
$\xib$-parameterized $\EM^\infty$-net $(\vphi_{1,\eps}(.,\xibR))_\eps$ also
satisfies conditions (i), (ii) in Lemma \ref{uni_lem}, as does
$\vphi_{2,\eps}(.,\xi)\cdot t^2 / r_\eps^2 :=  L_{\eps,\xibR}
\vphi_{1,\eps}(.,\xibR) - a_0^\eps \vphi_{1,\eps}(.,\xibR) /
P_{1,\eps}(.,\xibR)$. Another integration by parts in $I_1^\eps(\xi)$ gives
\begin{multline*}
    I_1^\eps(\xi) =
        \int e^{-ix\xi} \vphi_{2,\eps}(x,\xi) \cdot u_\eps(x) \, dx \\
        + (\frac{r_\eps}{t})^2
        \int e^{-ix\xi} \frac{\vphi_{1,\eps}(x,\xibR)}{P_{1,\eps}(x,\xib)}
            P_\eps u_\eps(x) \, dx =:
        I_2^\eps(\xi) + \frac{r_\eps^2}{t^2} h_2^\eps(\xib,\xi).
\end{multline*}
Again, $t\mapsto h_2^\eps(\xib,t \xib)$ is seen to be rapidly decreasing
uniformly in $\xib$ with $\eps$-power $-K$ (the same $K$ as above).

Successively, after $k$ steps we arrive at
\[
    (\vphi u_\eps)\FT{\ }(\xi) = I_k^\eps(\xi) +
                \sum_{j=1}^k (\frac{r_\eps}{t})^j h_j^\eps(\xib,t\xib)
\]
where $t \geq r_\eps$ and $h_j^\eps(\xib,t\xib) = O(t^{-k} \eps^{-K})$ when $0
< \eps < \eps_0$ uniformly in $\xib$ ($1 \leq j \leq k$), and
\[
    I_k^\eps(\xi) = (\frac{r_\eps}{t})^k \cdot \int e^{-ix\xi}
    \vphi_{k,\eps}(x,\xibR) u_\eps(x)\, dx
\]
with $\linf{\vphi_{k,\eps}(.,\xibR)} = O(\eps^{-1})$ uniformly in $\xib$.
Since $r_\eps^k = O(\eps^{-1})$ and $\sup |u_\eps(x)| =
O(\eps^{-M})$ ($\sup$ over $x\in\supp(\vphi)$) for some $M$ we finally find
that
\[
    |(\vphi u_\eps)\FT{\ }(\xi)| = O(t^{-k} \eps^{-M-2}) + O(t^{-k} \eps^{-K})
        = O(t^{-k} \eps^{-N})
\]
uniformly in $\xib$ when $t\geq \max(1,r_\eps)$ and $N := \max(M+2,K)$ with
arbitrary $k$.  \hfill $\square$

\begin{rem} 
(i) The simple examples $P_\eps = \diff{x} - i/\eps$ and $Q_\eps =
  \eps \diff{x} - i$, both admitting $u_\eps(x) = \exp(ix/\eps)$ as solution to
  the homogeneous equation, show that neither the slow scale condition on the
    coefficients nor the $W_{sc}$-ellipticity can be dropped in Theorem
    \ref{WF_thm}.

(ii) We may use this opportunity to point out that Theorem \ref{WF_thm}
  establishes corresponding claims made earlier in \cite[Examples 2.1 and 4.1]{HK:01}
    and in \cite[Theorem 23 (i)]{HdH:01} independently of the references
  cited therein (cf.\ Example \ref{Char_ex}(iii)).
\end{rem}

%% file: sec7b.tex
\section{Solvability of PDEs with coefficients in $\gC$}

We present a necessary and sufficient condition on the symbol $P(\xi) =
\sum_{|\al|\leq m} a_\al \xi^\al$ ($a_\al\in\gC$) for the corresponding PDE
to be solvable in $\Omega$ for arbitrary compactly supported right-hand
sides. More precisely, we investigate the property
\begin{equation}\label{solv}
    \forall f \in \G_c(\Om): \exists u \in \G(\Om): \qquad P(D) u = f.
\end{equation}

A sufficient condition for general solvability on $\Om = \R^n$ has been
given in \cite[Theorem 2.4]{NPS:98}, in the special case of the
existence of fundamental solutions already in \cite{NP:98}.
We restate it here in a slightly simplified form. (The
simplification is an immediate consequence of the characterization
of invertible generalized numbers.)
\begin{thm}[\cite{NPS:98}, Theorem 2.4]\label{NPSsolv}
Assume that there is $\xi_0\in\R^n$ such that
    \begin{equation}\label{NPScond}
        P_m(\xi_0) \text{ is invertible in } \gC.
    \end{equation}
Then $P(D)$ is solvable in $\R^n$, i.e.,
    \begin{equation}
        \forall f \in \G(\R^n): \exists u \in \G(\R^n): \qquad P(D) u = f.
    \end{equation}
\end{thm}
The proof in \cite{NPS:98} uses a complex Fourier integral representation
(within a partition of unity) of a solution candidate.
We will take a different approach while restricting to the case of
compactly supported right-hand side; on the other hand, we thereby gain a
relaxation of condition (\ref{NPScond}), which will, in addition, turn out to
be a characterization of solvable operators.

Before stating the new solvability condition we briefly discuss the relations
among various properties of the symbol. First, we note the simple fact that
condition
(\ref{NPScond}) is
implied by any of the ellipticity conditions introduced above.
\begin{prop} Every W-elliptic symbol $P(\xi)$ satisfies (\ref{NPScond}).
\end{prop}
\begin{proof} Recall from Proposition \ref{principalpartprop} that
W-ellipticity of
$P$ is equivalent to S-ellipticity of the principal part $P_m$. Hence there is
$N$, $R$, and $\eps_0$ such that $|P_{\eps,m}(\xi)| \geq \eps^N (1 + |\xi|)^m$
when $0 < \eps < \eps_0$ and $|\xi| \geq R$. Picking $\xi_0$ with $|\xi_0| \geq R$
arbitrary we obtain (\ref{NPScond}), after readjusting $N$ and $\eps_0$
accordingly.
\end{proof}
Clearly, condition (\ref{NPScond}) is strictly weaker than W-ellipticity, as
can be seen from the example $P(\xi_1,\xi_2) = \xi_1 - \xi_2$.

The following example shows that (\ref{NPScond}) is not necessary for
solvability.
\begin{ex}\label{zd_ex}
Let $P(\xi) = a \xi + i$, where $0 \not= a\in \gR$ is the zero divisor with
the following representative
\[
    a_\eps = \begin{cases}
            0 & \text{if $1/\eps \in \N$},\\
            1 & \text{otherwise}.
            \end{cases}
\]
Then $P_1(\xi) = a \xi$ cannot be invertible (for any $\xi$),
but $P(D) u = f \in\G$ 
is always solvable. A solution is given by the class of the representative
\[
    u_\eps(x) = \begin{cases}
            - i f_\eps(x) & \text{if $1/\eps \in \N$},\\
            i \int_0^x e^{(x-y)} f_\eps(y)\, dy & \text{otherwise}.
            \end{cases}
\]
Note, however, that here $|P(\xi)|^2 = a^2 \xi^2 + 1$ is invertible (for
arbitrary $\xi$) thanks to the lower order term.
\end{ex}

The key property of the (generalized) symbol $P(\xi)$, which will turn out to
be equivalent to (\ref{solv}), is specified in terms of its associated
(generalized) temperate weight function, which we define in analogy with
\cite[2.1, Example 2]{Hoermander:63} (see also \cite[Example
10.1.3]{Hoermander:V2})
\begin{defn} If $P$ is the symbol of a PDO of order $m$ with coefficients in
$\gC$ we define $\Pq := \sum_{\al\leq m} \d^\al P \cdot \overline{\d^\al P}
\in \G(\R^n)$ or, alternatively, in terms of representatives
\begin{equation}
  \label{PSchlange} \Pq_\eps(\xi) := \sum_{\al\leq m} |\d_\xi^\al
     P_\eps(\xi)|^2.
\end{equation}
\end{defn}
(Note that, contrary to the classical case, we avoid taking the square root,
but prefer to stick closely to the classical notation.)

\begin{lemma}\label{Pinv} If there is $\xi_0\in\R^n$ such that $\Pq(\xi_0)$
is invertible (in $\gC$) then $\Pq$ is invertible in $\G(\R^n)$ and $\Pt :=
\sqrt{\Pq}$ is a well-defined Colombeau function. More precisely, there
exist $d > 0$, $ N \geq 0$, $\eps_0 \in (0,1)$ such that
\begin{equation}\label{Pinv_lb}
    \Pq_\eps(\xi) \geq \eps^N (1 + d |\xi_0 - \xi|)^{-2m} \qquad
        \forall \xi\in\R^n, \forall \eps \in (0,\eps_0).
\end{equation}
\end{lemma}
\begin{proof} By \cite[Equation (2.1.10)]{Hoermander:63} there is a constant
$d > 0$ (independent of $\eps\in (0,1]$) such that
\[
    \Pq_\eps(\xi+\eta) \leq (1 + d |\eta|)^{2m} \Pq_\eps(\xi)
    \qquad \forall \xi, \eta \in \R^n, \forall \eps\in (0,1].
\]
Since $\Pq(\xi_0)$ is invertible we have 
for
some $N > 0$ and $\eps_0 \in (0,1]$ that $\eps^N \leq \Pq_\eps(\xi_0)$ when
$\eps \in (0,\eps_0)$. Therefore, substituting $\eta = \xi_0 - \xi$ in the
inequality above we obtain
\[
    \eps^N \leq \Pq_\eps(\xi_0) \leq (1 + d |\xi_0 - \xi|)^{2m} \Pq_\eps(\xi)
\]
for all $\xi\in\R^n$ and $\eps \in (0,\eps_0)$. This proves (\ref{Pinv_lb})
and shows that the square root of $\Pq_\eps$ is smooth (and moderate).
Moreover, (\ref{Pinv_lb}) yields the invertibility of $\Pq$ as a generalized
function on $\R^n$.
\end{proof}

The symbol in Example \ref{zd_ex} defines an invertible weight
function (e.g.\ $\Pq(0) = a^2 + 1$), whereas condition (\ref{NPScond}) is
not satisfied. The following proposition shows that, in general, the
invertibility of $\Pq$ is a (strictly) weaker condition.
\begin{prop} If $P(\xi)$ is a generalized symbol satisfying (\ref{NPScond})
then its associated weight function $\Pq$ is invertible.
\end{prop}
\begin{proof} Observe that
\[
    \frac{1}{n} \Big( \sum_{|\al|=m} |a^\eps_\al| \Big)^2 \leq
    \sum_{|\al|=m} (\al !)^2 |a^\eps_\al|^2 =
    \sum_{|\al|=m} |\d_\xi^\al P_\eps(\xi)|^2 \leq \Pq_\eps(\xi).
\]
Let $\xi_0\in \R^n$ be arbitrary. There is $C > 0$ (dependent only
on $n$ and $m$) such that
\[
    |P_{\eps,m}(\xi_0)| = |\sum_{|\al|=m} a_\al^\eps \xi_0^\al| \leq
        C |\xi_0|^m \sum_{|\al|=m} |a_\al^\eps|.
\]
Thus the invertibility of $P_m(\xi_0)$ implies the invertibility of
$\Pq_\eps(\xi_0)$.
\end{proof}

\begin{thm}\label{solv_thm} Assume that $\Pq$ is invertible at some
$\xi_0\in\R^n$. Then for every $f\in\G_c(\Om)$ there is a solution
$u\in\G(\Om)$ to the equation $P(D) u = f$.
\end{thm}
\begin{proof}  Recall the definition of the Banach spaces
$\B_{p,k}$ (cf.\ \cite[Ch.II]{Hoermander:63} or \cite[Ch.10]{Hoermander:V2}),
where $1 \leq p \leq \infty$ and $k$ is a temperate weight
function (cf.\ \cite[Ch.II]{Hoermander:63}). These are the spaces of
temperate distributions on $\R^n$ given by
\[
    \B_{p,k} = \{ v \in\S'(\R^n) \mid k\cdot\FT{u}\in \L^p(\R^n) \}
\]
and equipped with the norm $\| v \|_{p,k} = \|k\FT{u}\|_{L^p} / (2\pi)^{n/p}$.
Furthermore, the
corresponding local space $\B_{p,k}^{loc}$ consists of the distributions
$w\in\D'(\R^n)$ such that $\vphi w \in \B_{p,k}$ for every $\vphi \in
\D(\R^n)$.

We apply \cite[Theorem 3.1.1]{Hoermander:63} to obtain a regular fundamental
solution for the operator $P_\eps$, for each $\eps\in (0,1]$. More precisely,
we have the following. Let $c > 0$ be arbitrary; there is a constant $C > 0$,
depending only on $n$, $m$,
and $c$, such that $\forall \eps \in (0,1]$ $\exists E_\eps \in
\B_{\infty,\Pt_\eps}^{loc}$ such that $P_\eps(D) E_\eps = \de_0$ and
$E_\eps / \cosh(c|.|) \in \B_{\infty,\Pt_\eps}$ with norm estimate
\begin{equation}\label{E_norm_est}
    \| E_\eps / \cosh(c|.|) \|_{\infty,\Pt_\eps} \leq C.
\end{equation}

Choose a representative $(f_\eps)_\eps \in f$ and $K_0\Subset \Om$
such that $\supp(f_\eps) \subseteq K_0$ for all $\eps$. Define
\[
    u_\eps := (E_\eps * f_\eps)\vert_\Om \in \Cinf(\Om), \eps\in (0,1].
\]
By construction, $P_\eps(D) u_\eps = f_\eps$, so it remains to be shown that
$(u_\eps)_\eps$ is moderate.

Let $K \Subset \Om$, $\al\in \N_0^n$, $x\in K$ arbitrary. From the
definition of $u_\eps$ we obtain
\begin{multline*}
    |\d^\al u_\eps(x)| = |E_\eps * \d^\al f_\eps (x)| =
     |\dis{E_\eps}{\d^\al f_\eps(x - .)}| \\
    = |\dis{\frac{E_\eps}{\cosh(c|.|)}  }%
        {\big( (\d^\al f_\eps)
            \cdot\cosh(c |x-.|)\big)  (x - .)}| =
        |\big( e_\eps * g^\al_{\eps,x}\big) (x)|
\end{multline*}
where we have put $e_\eps:= E_\eps / \cosh(c|.|)$ and $g^\al_{\eps,x}(y) :=
\d^\al f_\eps(y) \cosh(c |x-y|)$.
Note that we have $e_\eps\in \B_{\infty,\Pt_\eps}$ and $g^\al_{\eps,x} \in
\D(K_0) \subset  (\B_{1,1 / \Pt_\eps} \cap \E'(K_0))$; now
\cite[Theorem 2.2.6]{Hoermander:63} yields  that  $e_\eps * g^\al_{\eps,x}
\in \B_{1,1} = \F^{-1}(\L^1) \subset \L^\infty$ and hence we have
\begin{multline*}
    |\big( e_\eps * g^\al_{\eps,x}\big) (x)|
        \leq \linf{ e_\eps * g^\al_{\eps,x}} \\
        \leq \| e_\eps * g^\al_{\eps,x}\|_{1,1} \leq
    \|e_\eps \|_{\infty,\Pt_\eps} \| g^\al_{\eps,x}\|_{1, 1 / \Pt_\eps}  \leq
        C \| g^\al_{\eps,x}\|_{1, 1 / \Pt_\eps}.
\end{multline*}
We have to establish a moderate upper bound for the last factor, that is
\begin{multline*}
    \| g^\al_{\eps,x}\|_{1, 1 / \Pt_\eps} =
        \lone{\FT{g^\al_{\eps,x}} / \Pt_\eps}
    = \int \frac{|\FT{g^\al_{\eps,x}}(\xi)|}{\Pt_\eps(\xi)} \, d\xi \\
    \leq \eps^{-N/2}
    \int |\FT{g^\al_{\eps,x}}(\xi)|\cdot (1 + d|\xi_0-\xi|)^m \, d\xi, \qquad
        \eps\in (0,\eps_0),
\end{multline*}
where we have made use of Lemma \ref{Pinv} (and the notation there).

A direct calculation, using Leibniz' rule, the support properties of $f_\eps$,
and noting that the factor $\cosh(c|x-.|)$ is $\eps$-independent, shows that
the family $(g^\al_{\eps,x})_{x\in K}$ has moderate upper
bounds (with respect to $\eps$) in every semi-norm of $\S(\R^n)$. By the
continuity of the Fourier transform we conclude that for every $l\in \N$ there
is $M \geq 0$ and $C_l > 0$ such that
\[
    |\FT{g^\al_{\eps,x}}(\xi)| \leq C_l (1 + |\xi|)^{-l} \eps^{-M},
\]
when $\eps$ is sufficiently small (and uniformly in $x\in K$).
Choosing $l > m + n$ we may use this bound
in the integrand above and arrive at
\[
    \| g^\al_{\eps,x}\|_{1, 1 / \Pt_\eps} \leq C' \eps^{-(2M + N)/2},
\]
where $C'$ is some constant depending only on $n$, $m$, $l$, and $\xi_0$, and
$\eps$ is small.

Combining these estimates, we have shown that
\[
    \sup_{x\in K} |\d^\al u_\eps(x)| \leq C C' \eps^{-(2M+N)/2}
\]
if $\eps$ is sufficiently small.
\end{proof}


The above solvability proof was based on a careful analysis and extension of
classical (distributional) constructions. Curiously enough, we can show
a converse implication by simple reasoning with properties of the ring of
generalized numbers.

\begin{thm}\label{nec_thm} If $f\in\G(\Om)$ such that $f(x_0)$ is invertible
(in $\gC$) for some $x_0\in\Om$ and $P(D) u = f$ is solvable in $G(\Om)$ then
$\Pq$ is invertible (in $\G(\Om)$).
\end{thm}
\begin{proof} Assume that $\Pq$ is not invertible. Then due to Lemma
\ref{Pinv}, in particular, $\Pq(0)$ cannot be invertible. Since $\d^\al
P_\eps(0) = \al! a_\al^\eps$ we have
\[
    \Pq_\eps(0) = \sum_{|\al|\leq m} (\al!)^2  |a_\al^\eps|^2.
\]
Since non-invertible generalized numbers are
zero divisors  we may choose a representative of
$\Pq(0)$, say $(b_\eps)_\eps$, which vanishes on a zero sequence of
$\eps$-values, i.e., $b_{\nu_k} = 0$ ($k\in\N$) for some sequence $(\nu_k)_k
\in (0,1]^\N$ with $\nu_k \to 0$ as $k\to\infty$.
We have for all $q$ that $|b_\eps - \Pq_\eps(0)| = O(\eps^q)$ ($\eps\to 0$).
Since all terms in the above sum representation of $\Pq_\eps(0)$ are
nonnegative, we deduce that $|a^{\nu_k}_\al| = O(\nu_k^q)$  for all $q$
($k\to\infty$). (In fact, we may choose representatives of $a_\al$, $|\al|\leq
m$, which vanish along the same zero sequence.)

Define the generalized number $c\in\gR$ by the representative
\[
    c_\eps = \begin{cases}
            1 & \text{if $\eps = \nu_k$ for some $k\in\N$},\\
            0 & \text{otherwise}.
            \end{cases}
\]
Clearly, $c \not= 0$ but $c \cdot a_\al = 0$ ($|\al|\leq m$) by construction.
If $f(x_0)$ is invertible and $P(D) u = f$ we arrive at the following
contradiction
\[
    0 \not= c\cdot f(x_0) = c \cdot P(D) u (x_0) =
        \sum_\al c a_\al D^\al u(x_0) = 0.
\]
\end{proof}

\begin{cor} The solvability property (\ref{solv}) for $P(D)$ with
        generalized constant coefficients holds if and only if $\Pq$
        is invertible (at some point in $\R^n$).
\end{cor}



%% file: ho.bbl
\newcommand{\SortNoop}[1]{}
\begin{thebibliography}{10}

\bibitem{Colombeau:84}
J.~F. Colombeau.
\newblock {\em New Generalized Functions and Multiplication of Distributions}.
\newblock North-Holland, Amsterdam, 1984.

\bibitem{Colombeau:85}
J.~F. Colombeau.
\newblock {\em Elementary Introduction to New Generalized Functions}.
\newblock North-Holland, 1985.

\bibitem{Garetto:02}
C.~Garetto.
\newblock Pseudo-differential operators in algebras of generalized functions
  and global hypoellipticity.
\newblock {\em Preprint, Quaderno N.22/2002, Universit\`{a} di Torino}, 2002.
\newblock To appear in Acta Appl. Math. 2003.

\bibitem{GKOS:01}
M.~Grosser, M.~Kunzinger, M.~Oberguggenberger, and R.~Steinbauer.
\newblock {\em Geometric theory of generalized functions}.
\newblock Kluwer, Dordrecht, 2001.

\bibitem{Hoermander:63}
L.~H{\"o}rmander.
\newblock {\em Linear Partial Differential Operators}.
\newblock Springer-Verlag, Berlin, 1963.

\bibitem{Hoermander:V2}
L.~H{\"o}rmander.
\newblock {\em The Analysis of Linear Partial Differential Operators},
  volume~II.
\newblock Springer-Verlag, 1983.

\bibitem{Hoermander:V1}
L.~H{\"o}rmander.
\newblock {\em The Analysis of Linear Partial Differential Operators},
  volume~I.
\newblock Springer-Verlag, second edition, 1990.

\bibitem{Hoermann:99}
G.~H{\"o}rmann.
\newblock Integration and microlocal analysis in {C}olombeau algebras.
\newblock {\em J. Math. Anal. Appl.}, 239:332--348, 1999.

\bibitem{HdH:01}
G.~H{\"o}rmann and M.~V. de~Hoop.
\newblock Microlocal analysis and global solutions of some hyperbolic equations
  with discontinuous coefficients.
\newblock {\em Acta Appl. Math.}, 67:173--224, 2001.

\bibitem{HK:01}
G.~H{\"o}rmann and M.~Kunzinger.
\newblock Microlocal analysis of basic operations in {C}olombeau algebras.
\newblock {\em J. Math. Anal. Appl.}, 261:254--270, 2001.

\bibitem{HOP:03}
G.~H{\"{o}}rmann, M.~Oberguggenberger, and S.~Pilipovi{\'{c}}.
\newblock Microlocal hypoellipticity of linear partial differential operators
  with generalized functions as coefficients.
\newblock {\em Preprint {\tt arXiv:math.AP/0303248}}, 2003.

\bibitem{LO:91}
F.~Lafon and M.~Oberguggenberger.
\newblock Generalized solutions to symmetric hyperbolic systems with
  discontinuous coefficients: the multidimensional case.
\newblock {\em J. Math. Anal. Appl.}, 160:93--106, 1991.

\bibitem{NP:98}
M.~Nedeljkov and S.~Pilipovi\'c.
\newblock Hypoelliptic differential operators with generalized constant
  coefficients.
\newblock {\em Proc. Edinb. Math. Soc., II. Ser.}, 41:47--60, 1998.

\bibitem{NPS:98}
M.~Nedeljkov, S.~Pilipovi{\'{c}}, and D.~Scarpal{\'{e}}zos.
\newblock {\em The Linear Theory of Colombeau Generalized Functions}.
\newblock Addison Wesley Longman, 1998.

\bibitem{O:89}
M.~Oberguggenberger.
\newblock Hyperbolic systems with discontinuous coefficients: generalized
  solutions and a transmission problem in acoustics.
\newblock {\em J. Math. Anal. Appl.}, 142:452--467, 1989.

\bibitem{O:92}
M.~Oberguggenberger.
\newblock {\em Multiplication of Distributions and Applications to Partial
  Differential Equations}.
\newblock Longman Scientific {\&} Technical, 1992.

\bibitem{Rosinger:90}
E.~E. Rosinger.
\newblock {\em Non-Linear Partial Differential Equations. An Algebraic View of
  Generalized Solutions}.
\newblock North-Holland, Amsterdam, 1990.

\bibitem{Schwartz:54}
L.~Schwartz.
\newblock Sur l'impossibilit{\'{e}} de la multiplication des distributions.
\newblock {\em C. R. Acad. Sci. Paris}, 239:847--848, 1954.

\end{thebibliography}
